\documentclass[onefignum,onetabnum]{siamart220329}

\usepackage{amssymb}
\usepackage{graphicx}
\usepackage{ulem}
\usepackage{multirow}
\usepackage{cite} % order citations [1-3,5]
\usepackage{bm}
\usepackage{xcolor}
\usepackage{algorithm}
\usepackage{pgfplots}
\usetikzlibrary{pgfplots.groupplots}
\usepackage{algpseudocode}
\usepackage{tkz-euclide}
\usepackage{tikz}
\usepackage{siunitx}
\usepackage{setspace}
\usepackage{todonotes}
\usepackage[mathscr]{euscript}
\usepackage{contour}
\pgfplotsset{compat=1.18}

\newsiamthm{defi}{Definition}
\DeclareMathAlphabet{\mathpzc}{OT1}{pzc}{m}{it}

\DeclareMathAlphabet{\mathbf}{OT1}{cmr}{bx}{it}

\newcommand{\norm}[1]{\left\lVert#1\right\rVert}

\newcommand{\T}{\mathcal{T}}

\newtheorem{example}[theorem]{Example}
\let\oldexample\example
\renewcommand{\example}{\oldexample\normalfont}

\newtheorem{remarksimple}[theorem]{Remark}
\let\oldremarksimple\remarksimple
\renewcommand{\remarksimple}{\oldremarksimple\normalfont}
\newenvironment{remark}{\begin{remarksimple}}{\hfill$\diamond$\end{remarksimple}}
\newtheorem{experiment}[theorem]{Experiment}
\let\oldexperiment\experiment
\renewcommand{\experiment}{\oldexperiment\normalfont}

\makeatletter 
\@mparswitchfalse%
\makeatother
\normalmarginpar
\title{Randomized sketched TT-GMRES for linear systems with tensor structure}
\author{Alberto Bucci\thanks{Department of Mathematics, University of Pisa, Italy, \texttt{alberto.bucci@phd.unipi.it}.
} \and Davide Palitta\thanks{Dipartimento di Matematica, (AM)$^2$, Alma Mater Studiorum - Università di Bologna, 40126 Bologna, Italy, \texttt{davide.palitta@unibo.it}} \and Leonardo Robol\thanks{Department of Mathematics, University of Pisa, Italy,
\texttt{leonardo.robol@unipi.it} }}
\date{\today}
\begin{document}

\renewcommand{\thefootnote}{\fnsymbol{footnote}}
\maketitle \pagestyle{myheadings} \thispagestyle{plain}
\markboth{A.\ BUCCI, D.\ PALITTA, L.\ ROBOL}{RANDOMIZED TT-SGMRES FOR TENSOR LINEAR SYSTEMS}

\begin{abstract}
In the last decade, tensors have shown their potential as valuable tools for various tasks in numerical linear algebra. While most of the research has been focusing on how to compress a given tensor in order to maintain information as well as reducing the storage demand for its allocation, the solution of linear tensor equations is a less explored venue. Even if many of the routines available in the literature are based on alternating minimization schemes (ALS), we pursue a different path and utilize Krylov methods instead. The use of Krylov methods in the tensor realm is not new. However, these routines often turn out to be rather expensive in terms of computational cost and ALS procedures are preferred in practice. We enhance Krylov methods for linear tensor equations with a panel of diverse randomization-based strategies which remarkably increase the efficiency of these solvers making them competitive with state-of-the-art ALS schemes. 
The up-to-date randomized approaches we employ range from sketched Krylov methods with incomplete orthogonalization and structured sketching transformations to streaming algorithms for tensor rounding. The promising performance of our new solver for linear tensor equations is demonstrated by many numerical results. 

\end{abstract}


\maketitle

%{ \color{blue}
%E' plausibile che l'ordine sotto non sia quello giusto, intanto 
%stiamo segnando tutte le cose di cui vogliamo parlare.
%}

\section{Introduction}
In the last decade linear tensor equations of the form 
\begin{equation}\label{eq:main}
    \mathcal{A}x=b,
\end{equation}
where $\mathcal{A}$ is an operator acting on 
$\mathbb R^{n_1\times\cdots\times n_d}$ and $x$, $b$ are tensors of appropriate dimensions,
have come up as very useful tools for describing the discrete problems stemming from a large setting of diverse applications. For instance, in, e.g., quantum chemistry~\cite{LubichQuantum,QuantumDynamics} and financial mathematics~\cite{Finance1,Finance2} high-order, possibly stochastic and parametric integral and partial differential equations (PDEs) need to be solved. 
 The discretization of these problems often leads to equations of the form~\eqref{eq:main}; see, e.g.,~\cite{Bachmayr_2023} and the references therein. Similarly, equation~\eqref{eq:main} can be used to model problems in imaging~\cite{tensor4Imaging} and deep neural networks~\cite{Tensor4NN} as well.

In spite of the large range of application settings where equation~\eqref{eq:main} can be met, only a handful of efficient solvers for its solution have been proposed in the literature. Most of them build upon (alternating) optimization schemes~\cite{ALS1,ALS2,ALS3} with AMEn~\cite{Amen}
and DMRG \cite{dmrg} being two of the most prominent representatives in this class of solvers. In~\cite{Multigrid} a multigrid procedure for~\eqref{eq:main} is proposed whereas in~\cite{TT-GMRES} a tensor-based implementation of the Generalized Minimal Residual (GMRES) method~\cite{GMRES} is presented and further studied in~\cite{TTGMRES2}. The numerical performance of some of these routines on multicore architectures has been recently investigated in~\cite{Comparison}.

In this paper, we assume that all the quantities in~\eqref{eq:main} are given in the \emph{tensor-train} (TT) format~\cite{TT_decomposition}. Indeed, this is one of the most suitable formats for representing (very) high-dimensional problems. Many of the procedures we are going to employ are tailored to this tensor format. However, the whole machinery we present here can be probably adapted to other formats as well.

The aim of this work is to significantly improve over the TT-GMRES method presented in \cite{TT-GMRES} by enhancing it with several 
randomization-based techniques developed in the last years in 
numerical linear algebra. TT-GMRES is a Tensor-Train 
formulation of the classic GMRES method. In particular, the basis vectors of the constructed Krylov subspace are represented in terms of TT-tensors, and TT-arithmetic is adopted throughout the iterative scheme. The computational cost of any operation involving TT-tensors depends linearly on the number of modes $d$ of the terms at hand but at least quadratically on their tensor rank; see~\cite[Section 4]{TT_decomposition}. Therefore, maintaining a small TT-rank during all the TT-GMRES iterations is crucial to obtain an affordable numerical scheme.
Unfortunately, both the application of the linear operator $\mathcal{A}$ in~\eqref{eq:main} and the orthogonalization step within TT-GMRES remarkably increase the TT-rank of the basis vectors. A low-rank truncation is thus performed after each of these steps to maintain the TT-ranks under control; see~\cite{TT-GMRES} and section~\ref{TTGMRES} for further details. As most 
Krylov methods in a low-rank (tensor) setting, the need 
to deal with repeated truncations can severely
affect the performance of the overall Krylov method; see, e.g., \cite{LowRankKrylov1,LowRankKrylov2} for details and analysis on some low-rank Krylov methods. 

We show that randomization can be a strong ally in 
this setting. First, we design a TT variant of the so-called 
sketched GMRES (sGMRES)~\cite{nakatsukasa2021fast}. This allows us to 
perform only a partial, incomplete reorthogonalization of the basis TT-vectors, with a consequent reduction in their TT-ranks, but still avoiding a drastic delay in the convergence of the underlying Krylov scheme. In addition to remarkably decreasing the overall computational efforts, the incomplete reorthogonalization step allows us to avoid storing the whole basis at all. While all the basis TT-vectors are clearly not necessary during the 
partial orthogonalization step, we show that their allocation can be avoided also to retrieve the final solution. In particular, we store and utilize only sketches of the basis vectors thanks to the employment of streaming low-rank approximation schemes \cite{kressner2023streaming,parallel_TT}. Notice that this is in contrast with different state-of-the-art Krylov-based procedures employing incomplete orthogonalization where the final solution is often retrieved by a so-called \emph{two-pass} strategy, namely a second Arnoldi step is performed at the end of the iterative procedure. 

All these different tools and ideas have a non-trivial interplay that we analyze in detail, especially from a computational point of view. 
We will show that our novel method is competitive and often more efficient than state-of-the-art linear solvers  for~\eqref{eq:main}. 
On the other hand, the many, diverse techniques we adopt make the derivation of sharp convergence bounds on the overall routine rather tricky and we thus leave this challenging, yet important, aspect to be studied elsewhere. 

Here a synopsis of the paper.
Section~\ref{sec: Background} sees some background material. In particular, we recall the general framework of sGMRES for (standard) linear systems, the TT-format, and TT-GMRES in section~\ref{sec: sGMRES},~\ref{tensor-trains}, and~\ref{TTGMRES}, respectively. The main contribution of this paper is illustrated in section~\ref{Sketched TT-GMRES} where we derive a sketched version of TT-GMRES (TT-sGMRES). All the randomization-based enhancements we equip TT-sGMRES with are presented in the following subsections. As any Krylov technique, when applied to poorly conditioned systems 
also our novel randomization-enhanced TT-sGMRES needs to be preconditioned to get a fast convergence in terms of number of iterations. This aspect is discussed in section~\ref{sec: preconditioning}. In section~\ref{sec: numerical experiments} a panel of diverse numerical results illustrates the potential of our procedure also when compared with different state-of-the-art techniques. The paper ends with some conclusions in section~\ref{sec: conclusions}.

\section{Background}
\label{sec: Background}
%----------------------------

In this section we provide a concise description of two essential ingredients for the construction of sketched TT-GMRES: the sketched GMRES method, and TT-GMRES, together with the main aspects of the TT-format. We only describe what is 
necessary for this paper, and we refer the reader
to \cite{inexact_gmres,simoncini_inexact} for further 
details on the former, and to \cite{TT-GMRES} for the latter. 

\subsection{Randomized sketching and GMRES}
\label{sec: sGMRES}

The Generalized Minimal Residual method (GMRES)~\cite{GMRES} is a classic iterative scheme for the numerical solution of large-scale, nonsymmetric systems of linear equations. Given a matrix $A\in\mathbb{R}^{n\times n}$ and a vector $b\in\mathbb{R}^n$ the algorithm approximates the solution to the linear system $Ax= b$. In particular, starting from an initial guess $x_0$, a solution $x_k$ of the form 
\begin{equation}\label{eq:GMRESsol}
 x_k=x_0+V_ky_k,   
\end{equation}
 is sought. The columns of the matrix $V_k=[v_1,\ldots,v_k]\in\mathbb{R}^{n\times k}$ form an orthonormal basis of the $k$-th Krylov subspace
\begin{equation}\label{eq:Krylov_def}
   \mathcal K_k(A, r_0)=\text{span}\{r_0,Ar_0,\ldots,A^{k-1}r_0\}, 
\end{equation}
where $r_0=b-Ax_0$ denotes the initial residual. The vector $y_k\in\mathbb{R}^k$ in~\eqref{eq:GMRESsol} solves the least squares problem
\begin{equation}\label{eq:minres}
y_k = \mathrm{argmin}_y \|AV_k y-r_0\|_2.
\end{equation}

If the basis $V_k$ is constructed by the \emph{full} Arnoldi method, namely an Arnoldi method where a full orthogonalization of the basis vectors is performed, the celebrated Arnoldi relation holds true, i.e.,
\begin{equation}\label{eq:ArnoldiRel}
    AV_k=V_{k+1}\underline{H}_k=V_kH_k+h_{k+1,k}v_{k+1}e_k^T,
\end{equation}
where $\underline{H}_{k}\in\mathbb{R}^{(k+1)\times k}$ collects the orthonormalization coefficients and $H_{k}\in\mathbb{R}^{k\times k}$ is its principal square submatrix;
see, e.g.,~\cite{saad}.

Thanks to orthogonality of $V_k$ the computation of $y_k$ in~\eqref{eq:minres} simplifies as 
\begin{equation}\label{eq:minres2}
y_k = \mathrm{argmin}_y \|AV_k y-r_0\|_2=\mathrm{argmin}_y \|\underline{H}_k y-\beta e_1\|_2, \quad \beta=\|r_0\|_2.
\end{equation}
Moreover, the current residual norm $\|Ax_k-b\|_2$ can be cheaply computed; see, e.g.,~\cite[Proposition 6.9]{saad}. GMRES terminates whenever $\|Ax_k-b\|_2$ satisfies a certain threshold condition. Otherwise, the Krylov subspace~\eqref{eq:Krylov_def} is expanded by computing a new basis vector and the scheme continues iteratively.

Many of the practical features and theoretical properties of GMRES 
depend on the orthogonality of the Krylov basis $V_k$. However, maintaining a fully orthogonal $V_k$ often 
becomes the bottleneck in practical computations, unless convergence is fast. 

Several strategies have been proposed over the years to mitigate this issue. A 
standard approach is to \emph{restart} either explicitly~\cite[Section 6.5.6]{saad} or implicitly by deflated restarting \cite{deflated-morgan}. 
Another option to lower the computational cost of the orthogonalization step is to perform an \emph{incomplete} orthogonalization, namely the new basis vector $v_k$ is explicitly orthogonalized only with respect to a certain number $\ell$ of previously computed $v_i$s; see, e.g.,~\cite[Section 6.5.7]{saad}. A strategy with a different flavor is  \emph{preconditioning}, where the original problem is implicitly transformed into a problem for which GMRES converges in fewer iterations. Reducing the number of iterations clearly lowers the cost of the orthogonalization as well. However, selecting the right preconditioner may be tricky, problem-dependent, and its application time consuming. 
While these approaches all share similar goals, they are often applied independently of each other. In the following, we will show that for tensor equations of the form~\eqref{eq:main} it is often sensible to integrate the aforementioned techniques to attain a very efficient solution scheme.

% \begin{algorithm}
%\caption{GMRES} 
%\label{alg:GMRES}
%\vspace{2mm}
%\textbf{Input} Matrix $A\in\mathbb{R}^{n\times n}$, right hand side $b\in\mathbb{R}^n$, initial guess $x_0\in\mathbb{R}^n$, $d$ maximum basis dimension, tolerance $\mathrm{tol}$.\\
%\textbf{Output} Approximate solution $\hat{x}$  such that $\|A \hat{x} - b\|_2 \leq \mathrm{tol}$.\\

%\rr{TODO}
%\vspace{2mm}
%\end{algorithm}

% { \color{red} 
% \begin{itemize}
%     \item Orthgonalization is often the dominant cost. 
%     \item This can be mitigated with restarts, or with randomization. 
% \end{itemize}}

 At this point we focus on the incomplete orthogonalization GMRES scheme. For this GMRES variant, the basis $V_k$ is no longer orthogonal. However, the Arnoldi relation~\eqref{eq:ArnoldiRel} still holds and the vector $y_k$ may still be computed as 
\begin{equation}\label{eq:argmin2}
    y_k =\mathrm{argmin}_y \|\underline{H}_k y-\beta e_1\|_2.
\end{equation}

 Nevertheless, due the nonorthogonality of the basis, $y_k \neq \mathrm{argmin}_y \|AV_k y-r_0\|_2.$
 It is well-known that this drawback often leads to a delay in the convergence of the solution scheme, in general. 
 However, in the recent literature, it has been shown that when combined with sketching techniques, GMRES with incomplete orthogonalization is often able to retrieve the rate of convergence of the fully orthogonal procedure;
see~\cite{nakatsukasa2021fast}. 

The integration of sketching and GMRES with incomplete orthogonalization, called \emph{sketched} GMRES (sGMRES), makes use of oblivious subspace embeddings as sketching matrices. In particular, given a $k$-dimensional subspace $\mathcal{V}_k$, a linear transformation $S\in \mathbb{R}^{s\times n}$, with $s>k$, is a subspace embedding with distortion $\varepsilon\in[0,1)$ for $\mathcal{V}_k$
if, for any $v\in \mathcal V_k$, we have
\begin{equation} \label{eq:subspace_embedding}
    (1-\varepsilon)\|v\|_2^2 \leq \|Sv\|_2^2 \leq (1+\varepsilon) \|v\|_2^2;
\end{equation}
see, e.g.,~\cite{sarlos2006improved,DrineasMahoneyMuthukrishnan2006,woodruff2014sketching}. Notice that the sketching matrix induces the semidefinite inner product $x^TS^TSy$. It can be shown that this is indeed an actual inner product on the space $\mathcal{V}_k$ for which $S$ is an $\varepsilon$-subspace embedding; see, e.g.,~\cite[Section 3.1]{BalabanovNouy2019}.

In our case, the space $\mathcal{V}_k$ corresponds to the Krylov subspace~\eqref{eq:Krylov_def} which is clearly not known a priori. Therefore, we will need to employ \emph{oblivious} subspace embeddings (OSEs) in our work. These are particular transformations $S$ that can be constructed by solely knowing the dimension of the subspace to be embedded and such that~\eqref{eq:subspace_embedding} holds with high probability. 
Common choices for oblivious subspace embeddings are, e.g.,  Gaussians, for their theoretical guarantees, or subsampled trigonometric transforms since they allow for fast application; see, e.g.,~\cite{randomproof}.

In~\cite{nakatsukasa2021fast}, the authors integrate sketching and GMRES by replacing the selection of $y_k$ in~\eqref{eq:minres} by the following condition
\begin{equation}\label{eq:sketched_argmin}
    y_k = \mathrm{argmin}_y \|SAV_k y - Sr_0\|_2,
\end{equation}
where the basis $V_k$ of the Krylov subspace is computed by an Arnoldi scheme with incomplete orthogonalization. Due to the lack of an Arnoldi-like relation for the sketched quantities in~\eqref{eq:sketched_argmin}, in~\cite{nakatsukasa2021fast} the authors compute $y_k$ by performing 
\begin{equation}\label{eq:pseudo_inv}
    y_k=(SAV_k)^\dagger Sr_0. 
\end{equation}
In Algorithm~\ref{alg:sGMRES} we report the overall sGMRES algorithm. 
\begin{algorithm}[t]\caption{sGMRES}
\begin{algorithmic}[1]
%\setstretch{1.2}
\smallskip

\Statex \textbf{Input:} Matrix $A\in\mathbb{R}^{n\times n}$, right-hand side $b\in\mathbb{R}^n$, initial guess $x_0\in\mathbb{R}^n$, maximum basis dimension $\texttt{maxit}$, sketching $S\in\mathbb{R}^{s\times n}$, incomplete orthogonalization parameter $\ell$, tolerance $\texttt{tol}$.  
\Statex \textbf{Output:} Approximate solution $x_k$  such that $\|S(Ax_k-b)\|\leq \| Sb \| \cdot \texttt{tol}$
\smallskip

\State Set $r_0=b-Ax_0$, $V_1=v_1=r_0/\|r_0\|$, $W_0 = []$
 \For{$k=1,\ldots,\texttt{maxit}$}
\State Compute $\widetilde v=Av_k$
\State Update $W_k=[W_{k-1},S\widetilde v]$
\For{$i=\max\{1,k-\ell+1\},\ldots,k$}
\State Set $\widetilde v=\widetilde v- v_ih_{i,k}$, where $h_{i,k}= \widetilde v^Tv_i$ 
\EndFor
\State Set $h_{k+1,k}=\|\widetilde v\|$ and $v_{k+1}=\widetilde v/h_{k+1,k}$
%\State Update skinny QR: $Q_{k+1}T_{k+1}=[SV_{k},Sv_{k+1}]$ \label{white:line1}
%\State Update $\widehat H_k=T_kH_kT_k^{-1}$ and set $\widehat h=(h_{k+1,k}/\tau_k)t_{k+1}$ \label{white:line2}
\State Compute $y_k$ as the solution to~\eqref{eq:pseudo_inv}
%\State Compute $\|r_k\|$ as in~\eqref{eq:sketched_residualnorm}
\If{$\|W_ky_k-Sr_0\|\leq \|Sb\| \cdot \texttt{tol}$}
\State Go to line~\ref{alg:final_line}
\EndIf
\State Set $V_{k+1}=[V_k,v_{k+1}]$
\EndFor
\State Set $x_k=x_0+V_ky_k$\label{alg:final_line}
\end{algorithmic}\label{alg:sGMRES}
\end{algorithm}

\subsection{Tensor-Train decomposition}\label{tensor-trains}

A tensor $\mathcal{T}$ of size $n_1 \times n_2 \times \dots \times n_d$ is in the TT-format if it can be written element-wise as

\begin{equation}\label{eq:tt-format}
\mathcal{T}[i_1,\dots, i_d] = \sum_{\ell_1=1}^{r_1}\dots \sum_{\ell_{d-1}=1}^{r_{d-1}} C_1[1, i_1, \ell_1] C_2[\ell_1, i_2, \ell_2]\dots C_d[\ell_{d-1}, i_d, 1].
\end{equation}

The third-order tensors $C_{\mu}$ of size $r_{\mu-1}\times n_\mu \times r_\mu$ are the TT-cores (where $r_0 = r_d = 1$). Using MATLAB notation, relation \eqref{eq:tt-format} can be written compactly as a product of $d$ matrices (where the
first and last matrices collapse to a row and column vector, respectively) as follows:
\[
\mathcal{T}[i_1,\dots, i_d] = C_1[1, i_1, :] C_2[:, i_2, :] \dots C_d[:, i_d, 1].
\]

% For a vector $v\in\mathbb{R}^n$, where $n=n_1\cdots n_d$, it is sufficient to map its index $i$ to the tuple $(i_1,\dots, i_d)$ and then approximate it in tensor train format.
% Possibile proposta per spezzettare questa sezione in modo 
% più digeribile, e più facile da ricercare a posteriori. 
In order to establish the notation, we 
briefly recall the basic operations on tensors 
that will be used in the next sections. \medskip 

\noindent \textbf{Unfoldings.}
The unfolding $\T_{\leq \mu}$ is one of the many ways to matricize a tensor; it is a matrix
of size $\prod_{k=1}^{\mu} n_k \times  \prod_{k = \mu+1}^d n_k$ obtained from merging the first $\mu$  modes of $\T$  into
row indices and the last $d-\mu$  modes into column indices.

\medskip 
\noindent \textbf{Interface matrices.}
Each unfolding can be factorized in a low-rank way as $C_{\leq \mu}C_{> \mu}^T$ where
\[
C_{\leq_\mu}  \in  \mathbb{R}^{(n_1 \cdots n_\mu) \times  r_\mu}\quad   \text{and} \quad  C_{>\mu}  \in  \mathbb{R}^{(n_{\mu +1}\cdots n_d)\times r_\mu}.
\]
These are sometimes called interface matrices. \medskip 

The tuple $(r_1,\dots, r_{d-1})$ is called the TT-representation rank of the Tensor-Train defined in \eqref{eq:tt-format} and it determines the
complexity of working with a TT. For instance, storing a tensor in TT-format requires storing the $O(dnr^2)$ entries of its TT-cores, where $n:= \max_\mu(n_\mu)$ and $r \approx r_\mu$ for all $\mu = 1, \ldots, d$.\footnote{For the sake of readability, we will often 
make the simplifying assumption that all TT-ranks can be 
estimated by a single scalar $r$, and the 
dimensions $n_\mu$ by $n_\mu \approx n$. This will make writing 
computational complexities much easier. The most general 
result can usually be recovered by replacing terms 
such as $dr^j$ with $\sum_{\mu = 1}^d r_\mu^j$, and analogously for 
the $n_\mu$. }
Any tensor can be trivially written in the TT-format by choosing the TT-representation ranks sufficiently large.
The TT-representation rank of a particular tensor $\mathcal{T}$ is by no means unique but there exists an (entry-wise) minimal value which is called the TT-rank of $\mathcal{T}$. The minimal value for $r_\mu$ equals the matrix rank of $\mathcal T_\mu$.
In the rest of the paper will not distinguish between TT-rank and TT-representation rank and simply call
$(r_1,\dots, r_{d-1})$ the TT-rank of the tensor $\mathcal{T}$ once relation \eqref{eq:tt-format} is satisfied for some cores $C_\mu$.

When dealing with vectors in tensor-train format, to simplify the matrix-vector products, it is preferable to write matrices in the tensor-train operator format.

A matrix $A$ of size $m\times n = (m_1\times \dots \times m_d)\times (n_1 \times \dots \times n_d)$ is in the operator TT-format if it can be written element-wise as
\begin{equation}\label{eq:operator tt-format}
    A[i_1,\dots,i_d, j_1,\dots, j_d] = D_1[1, i_1, j_1, :] D_2[:, i_2, j_2, :]\dots D_d[:, i_d, j_d, 1].
\end{equation}
Then, given a vector $v$ in TT-format with cores $C_k$s, 
to compute the cores $G_1, \dots, G_d$ of $y=Av$
it is possible to act on each core separately. In formulas
\[
G_k[(\ell_{k-1},\alpha_{k-1}),i_k, (\ell_{k},\alpha_{k})] = \sum_{j_k=1} D_k[\alpha_{k-1}, i_k, j_k, \alpha_k] C_k[\ell_{k-1}, j_k,\ell_k] .
\]
As we can see, the TT-ranks of the MatVec are bounded by the product 
of the TT-ranks of the matrix and the vector. In iterative schemes like GMRES, several applications of $\mathcal A$ are required; 
without rounding, this unavoidably leads to the TT-ranks becoming too large. Hence, a tensor rounding procedure, or compression, is needed.
A given tensor $\mathcal T$ is approximated by another tensor $\widetilde{\mathcal T}$ with minimal possible TT-ranks $(r_1, \dots, r_{d-1})$ with a
prescribed accuracy $\varepsilon$ (or a fixed maximal TT-rank $R$) if:
\[
 \|\mathcal T-\widetilde{\mathcal T}\|_F\leq \varepsilon \|\mathcal T\|_F \quad (\text{or}\hspace{1.5mm} r_{k}\leq R) .
\]
A quasi-optimal $\widetilde{\mathcal T}$ can be obtained by the TT-SVD algorithm
\cite{TT_decomposition} with $\mathcal O(dnr^3)$ complexity.
This is based on performing QR decomposition and truncated 
SVDs of the interface matrices, exploiting the low-rank structure. 
Cheaper (and at the same time slightly less accurate) alternatives 
are available \cite{rounding_TT, kressner2023streaming, TT_decomposition,parallel_TT}, and are often 
based on randomization. 

In this work, we will focus 
on \emph{streamable} and \emph{randomized} 
rounding schemes, i.e., algorithms that 
allow us to find a low-rank representation of a sum 
of tensors 
$\mathcal T^{(1)} + \ldots + \mathcal T^{(m)}$ by performing 
preliminary contractions of the tensors $\mathcal T^{(k)}$, and 
reconstructing the low-rank approximation of their sum at a later 
stage. This choice will bring benefits in both speed and accuracy, 
and will be discussed in further detail in section~\ref{sec:final-solution}.

\subsection{TT-GMRES}\label{TTGMRES}

TT-GMRES \cite{TT-GMRES} is an extension of GMRES aimed at solving tensor equations of the form~\eqref{eq:main} in tensor-train format. The main distinction from the classic GMRES is in the representation of the basis ``vectors'' $v_k$s which are now given as tensor-train vectors. Moreover, TT-GMRES sees the incorporation of rounding steps throughout the process to maintain the TT-ranks of the $v_k$s within a specified threshold.

In \cite{TT-GMRES} a truncation strategy
based on the theory of inexact GMRES~\cite{simoncini_inexact} is suggested. Heuristically, employing this procedure often keeps the TT-ranks under control. However, there is no clear theoretical link between this strategy and the growth of the ranks. 
Further exploration and insights in this direction would undoubtedly yield valuable contributions to the field. Similarly, the truncations taking place after the Gram-Schmidt cycle can potentially destroy the orthogonality of the basis making the analysis even trickier. This issue has been studied in~\cite{LowRankKrylov1} in the case of low-rank Krylov methods for multiterm matrix equations.

In Algorithm~\ref{alg:TT-GMRES} we report the overall TT-GMRES scheme. In lines~\ref{TT-GMRES_round1} and \ref{TT-GMRES_round2}, 
$\Call{Round}{\mathcal{T}, \theta}$ denotes the TT-SVD from 
\cite{TT_decomposition} that performs a $\theta$-accurate low-rank truncation of the tensor $\mathcal{T}$. 

\begin{algorithm}[t]\caption{TT-GMRES}
\begin{algorithmic}[1]
%\setstretch{1.2}
\smallskip

\Statex \textbf{Input:} Tensor $\mathcal A\in\mathbb{R}^{n_1\times\ldots\times n_d}$, right-hand side $b$, initial guess $x_0$ in TT-format, maximum basis dimension $\texttt{maxit}$, tolerance $\texttt{tol}$.
\Statex \textbf{Output:} Approximate solution $x_k$  such that $\|\mathcal Ax_k-b\|\leq \| b \| \cdot \texttt{tol}$
\smallskip

\State Set $r_0=b-\mathcal Ax_0$, $\beta=\|r_0\|$, $V_1=v_1=r_0/\beta$ 
 \For{$k=1,\ldots,\texttt{maxit}$}
\State Set $\eta_k = \texttt{tol} / \norm{r_{k-1}}$
\State Compute $\widetilde v=\Call{Round}{\mathcal A v_k, \eta_k \cdot \texttt{tol}}$\label{TT-GMRES_round1}
\For{$i=1,\ldots,k$}
\State Set $\widetilde v=\Call{Round}{\widetilde v- v_ih_{i,k}, \eta_k \cdot \texttt{tol}}$, where $h_{i,k}= \widetilde v^Tv_i$ 
\label{TT-GMRES_round2}
\EndFor
%\State Set $\widetilde v=\Call{Round}{\widetilde v, \eta_k\cdot \mathrm{tol}}$\label{TT-GMRES_round2}
\State Set $h_{k+1,k}=\|\widetilde v\|$ and $v_{k+1}=\widetilde v/h_{k+1,k}$
\State Compute $y_k$ as the solution to~\eqref{eq:argmin2}
\State Compute $\|r_k\|=\|\underline{H}_ky_k-\beta e_1\|$ 
\If{$\|r_k\|\leq \| b \| \cdot\texttt{tol}$}
\State Go to line~\ref{alg:final_line2}
\EndIf
\State Set $V_{k+1}=[V_k,v_{k+1}]$
\EndFor
\State Set $x_k=x_0$
\For{$i = 1, \ldots, k$}
    \State $x_k = \Call{Round}{x_k + v_i \cdot (e_i^T y_k), \texttt{tol}}$
    \label{alg:final_line2}
\EndFor
\end{algorithmic}\label{alg:TT-GMRES}
\end{algorithm}

Thanks to the theory of inexact Arnoldi \cite{simoncini_inexact}, 
the roundings
at line~\ref{TT-GMRES_round1} and \ref{TT-GMRES_round2}
can be made more aggressive as the method converges, which 
helps to maintain the basis vectors of moderate ranks.
Nevertheless, the full orthogonalization step makes the overall procedure extremely time-consuming, in general. This is one of the reasons why TT-GMRES is not commonly employed for the solution of~\eqref{eq:main} and ALS procedures are often preferred. In the following sections we propose a sketched variant of TT-GMRES which, when equipped with a series of other randomization-based tools, turns out to be competitive with respect to state-of-the-art ALS schemes; see section~\ref{sec: numerical experiments}.

\section{Sketched TT-GMRES}\label{Sketched TT-GMRES}

The previous sections provided the necessary tools and theoretical background to facilitate the understanding of the sketched Tensor-Train GMRES (TT-sGMRES) method, which we present here.

The structure of the section is as follows. In Algorithm~\ref{alg:sTT-GMRES}, we begin by outlining the pseudocode for adapting the sGMRES algorithm to the TT-format, akin to the TT-GMRES approach given in Algorithm \ref{alg:TT-GMRES}.
The algorithm fundamentally expands upon sGMRES \cite{nakatsukasa2021fast}, adapting it to the TT-format similarly to how TT-GMRES in \cite{TT-GMRES} builds upon the GMRES method. This simple
generalization is not competitive with state-of-the art 
methods; hence, we delve into a series of refinements and techniques for its efficient implementation, that will turn 
it into a practical algorithm. 
In particular, we propose different techniques that exploit randomization to reduce the growth of the ranks, the memory requirements, 
and the cost of reorthogonalization; these techniques also reduce 
the cost and improve the stability of forming the final solution.  
Algorithm~\ref{alg:enhanced sTT-GMRES} summarizes these refinements in a detailed implementation.

% The sTT-GMRES algorithm can be seen as a traditional sGMRES algorithm implemented in TT-format, with additional randomized techniques to enhance its efficiency.
% We summarize these in Algorithm \ref{alg:sTT-GMRES}.

\begin{algorithm}[t]\caption{TT-GMRES}
\begin{algorithmic}[1]
%\setstretch{1.2}
\smallskip

\Statex \textbf{Input:} Tensor $\mathcal A\in\mathbb{R}^{n_1\times\ldots\times n_d}$, right-hand side $b$, initial guess $x_0$ in TT-format, maximum basis dimension $\texttt{maxit}$ , tolerance $\texttt{tol}$, sketching $S$, incomplete orthogonalization parameter $\ell$.  
\Statex \textbf{Output:} Approximate solution $x_k$  such that $\|\mathcal Ax_k-b\|\leq \| Sb \| \cdot \texttt{tol}$
\smallskip

\State Set $r_0=b-\mathcal Ax_0$, $\beta=\|r_0\|$ $V_1=v_1=r_0/\beta$, 
$W_0 = []$
 \For{$k=1,\ldots,\texttt{maxit}$}
\State Compute $\widetilde v=\Call{Round}{\mathcal A v_k, \nu_k\cdot \texttt{tol}}$\label{TT-GMRES_round12} \Comment{Adjustment for the tolerance, see \ref{sec:truncation}}
\State Update $W_k=[W_{k-1},S\widetilde v]$ \label{TT-sGMRES_updatesketching}
\For{$i=\max\{{1,k-\ell+1,\ldots,k}\}$}
\State Set $\widetilde v=\widetilde v- v_ih_{i,k}$, where $h_{i,k}= \widetilde v^Tv_i$  \label{TT-GMRES_orth}             \Comment{see~\ref{Randomized approximation of sums in TT-format} and  \ref{sec:reorthogonalization}}

\EndFor
\State Set $\widetilde v=\Call{Round}{\widetilde v, \eta_k\cdot \texttt{tol}}$\label{TT-GMRES_round22} \Comment{Adjustment for the tolerance, see \ref{sec:truncation}}

\State Set $h_{k+1,k}=\|\widetilde v\|$ and $v_{k+1}=\widetilde v/h_{k+1,k}$
%\State Update skinny QR: $Q_{k+1}T_{k+1}=[SV_{k},Sv_{k+1}]$ \label{white:line12}
%\State Update $\widehat H_k=T_kH_kT_k^{-1}$ and set $\widehat h=(h_{k+1,k}/\tau_k)t_{k+1}$ \label{white:line22}
\State Compute $y_k$ as the solution to~\eqref{eq:pseudo_inv}\label{sTT-GMRES_compute_y}
%\State Compute $\|r_k\|$ as in~\eqref{eq:sketched_residualnorm}
\If{$\|W_ky_k-Sr_0\|\leq \| Sb \| \cdot\texttt{tol}$}\label{sTT-GMRES_compute_resnorm}
\State Go to line~\ref{alg:final_line22}
\EndIf
\State Set $V_{k+1}=[V_k,v_{k+1}]$ 
\EndFor
\State Set $x_k = x_0$
\For{$i = 1, \ldots, k$}
    \State $x_k = \Call{Round}{x_k + v_i \cdot (e_i^T y_k), \texttt{tol}}$ \Comment{see \ref{sec:final-solution}}
    \label{alg:final_line22}
\EndFor
%\State Set $x_k=\Call{Round}{x_0+V_ky_k,\mathrm{tol}}$           \Comment{\rr{A: remove and change formula for sum}Sum computed by sketching, see \ref{sec:final-solution}}

\end{algorithmic}    \caption{Sketched TT-GMRES (TT-sGMRES) -- vanilla version}
    \label{alg:sTT-GMRES}

\end{algorithm}

%\begin{algorithm}[t]
%    \begin{algorithmic}[1]
%    \Procedure{sTTGMRES}{$A,b,x_0,S, k, \tau$}
%        \State $\eta \gets 10^{-1}$ \Comment{Adjustment for the tolerance, see \ref{sec:truncation}}
%        \State $r_0 \gets b - A x_0$
%        \State $v_1 \gets r_0 / \| r_0 \|_2$
%        \State $r_{0,S} \gets S r_0$ 
%           \Comment{Sketch of the initial residual}
%        \State $\mathrm{res} \gets \norm{r_{0,S}}$
%        \While{$\mathrm{res} > \tau $}
%            \State $w_i \gets \Call{Round}{A v_i, \eta \tau}$
%            \State $w_{i,S} \gets Sw_i$
%            \State $w_i \gets w_i - 
%               \sum_{j = 1}^{\min(k, i)}  \langle 
%                 v_{i-j+1}^T, w_i \rangle_S\ 
%                 v_{i-j+1}$
%              \Comment{see \ref{sec:reorthogonalization}}
%            \State $v_{i+1} \gets 
%              \Call{Round}{w_i / \norm{w_i}_2, \eta \tau}$
%            \State $W \gets \begin{bmatrix}
%                w_{1,S} & \dots & w_{i,S}
%            \end{bmatrix}$
%            \State $y \gets W^\dagger r_S$
%            \State $\mathrm{res} \gets \| r_{0,S} - W y \|$
%        \EndWhile
%        \State \Return 
%          \Call{Round}{$y_1 v_1 + \ldots + y_i v_i, \tau$}
%          \Comment{Sum computed by sketching, see \ref{sec:final-solution}}
%        \EndProcedure
%    \end{algorithmic}
%    \caption{Sketched TT-GMRES algorithm for the 
%  linear system $Ax = b$, with tolerance $\tau$, 
%  starting guess $x_0$, 
%  sketching $S$, and orthogonalization parameter $k$.}
%    \label{alg:sTT-GMRES2}
%\end{algorithm}

\subsection{Choice of structured sketchings}

The first aspect we discuss is the choice of the sketching $S\in \mathbb{R}^{s\times \prod_{k=1}^d n_k}$.  Notice that this transformation maps vectors in TT-format into standard vectors of $\mathbb{R}^s$. Therefore, no operations with sketched quantities as, e.g., the computation of $y_k$ in line~\ref{sTT-GMRES_compute_y}, involve tensor arithmetic.

Due to the huge number of columns, using  Gaussian transformations or subsampled trigonometric transforms for $S$ is prohibitively expensive and highlights the need for structure embeddings that exploit the TT structure of the vectors.

There are two natural ways to sketch a vector in TT-format, one based on the Kronecker product of matrices, and the other based on the Khatri-Rao product.
In particular, given a set of matrices  $S_1, \dots S_d$, with $S_k\in \mathbb{R}^{s_k\times n_k}$, and a TT-vector $\mathcal{T}$ with core tensors $C_k\in \mathbb{R}^{r_k\times n_k \times r_{k+1}}$, if we define $S_\otimes := S_1 \otimes \ldots \otimes S_d$, then the product $S_\otimes \mathcal{T}$ can be easily computed as it results in a TT-vector with cores $D_k = C_k \times_2 S_k$. In other words, the product is distributed across the cores, providing an exponential speed-up in the computation.
Notice that the transformation $S_\otimes$ maps vectors of length $\prod_{i=1}^dn_i$ into vectors of length $s = \prod_{i=1}^ds_i$. 
%Moreover, it can be shown that for $S_\otimes$ to be a $k$-dimensional oblivious subspace embedding. {\color{red} DP: qualcosa non va con questa frase!}

A different option is to draw matrices $S_k$ with the same number of rows and to opt for $S_\odot = S_1 \odot \dots \odot S_d$ where $\odot$ denotes the row-wise Khatri-Rao product, i.e., the $j$-th row of $S_\odot$ is the Kronecker product of the $j$-th rows of the matrices $S_k$s.
The advantage of this second operator is that its application on a TT-vector still splits across the cores, reducing the embedding cost; this computational gain comes at a minimal cost in embedding power \cite{Kolda}. For this reason in our algorithms we opt for the Khatri-Rao sketchings.
Regarding the distribution of the embeddings \( S_k \), we choose Gaussian embeddings to strengthen the theoretical guarantees. Specifically, each \( S_k \) is a Gaussian matrix with i.i.d. entries following \(\mathcal{N}(0, s^{-1/d})\) for appropriate scaling.

The selection of $s$ will be discussed in detail in 
section~\ref{sec:alltogether}. 

%For \( s \), we set it twice the maximum number of iterations of the algorithm. {\color{red}DP: I believe this last bit can go to Section~\ref{sec: numerical experiments}.}

\subsection{Truncation policy}
\label{sec:truncation}
One of the aspects that plays an important role in making Algorithm~\ref{alg:sTT-GMRES} competitive is the selection of the truncation tolerance for the rounding steps. Indeed, this must be able to avoid an excessive growth of the TT-ranks.

Algorithm~\ref{alg:sTT-GMRES} sees two main sources of rank growth: the application of $\mathcal A$ in line~\ref{TT-GMRES_round12} and the linear combinations of the basis vectors which occur both in the orthogonalization phase (line~\ref{TT-GMRES_orth}) and in the construction of the final solution (line~\ref{alg:final_line22}).
In \cite{TT-GMRES} the author suggests truncating the resulting tensors using the TT-SVD after each of these operations.
In particular, as noted in \cite{TT-GMRES}, the truncation taking place right after the matrix-vector product $\mathcal Av_k$ can be interpreted as an \emph{inexact} application of $\mathcal A$ to $v_k$.
Therefore, in principle, the theory of inexact Krylov methods can be employed to select suitable truncation parameters which do not jeopardize the convergence of the overall scheme.
The inexact GMRES method has been thoroughly examined by Szyld and Simoncini \cite{simoncini_inexact}, who introduce a progressively relaxed truncation policy. They prove that the accuracy in the application of $\mathcal A$ can be decreased gradually during the iterations. In particular,  if $\sigma_{\min}(\mathcal{A})$ denotes the smallest singular value of $\mathcal{A}$, then in \cite{simoncini_inexact} the authors suggest employing an iteration-dependent tolerance of the form
\begin{equation}\label{eq:eta_selection}
   \nu_k=\frac{\sigma_{\min}(\mathcal{A})}{\mathrm{maxit}\cdot\|r_{k-1}\|}. 
\end{equation}
In \cite{TT-GMRES} a similar value for the truncation in the rounding procedure is chosen. 

Notice that decreasing the accuracy in the application of $\mathcal{A}$ is equivalent to performing more aggressive low-rank truncations in our context. This is a rather crucial point as the TT-rank of the basis vectors $v_k$ increases with $k$ and being able to significantly reduce it in later iterations is thus extremely beneficial.

The proofs in  \cite{simoncini_inexact} strongly rely on the orthogonality of the basis $V_k$. However, the truncation taking place after the Gram-Schmidt step (line~\ref{TT-GMRES_orth} in Algorithm~\ref{alg:sTT-GMRES}) may potentially destroy the orthogonality of the basis, also in case of a full orthogonalization. This drawback should not get overlooked in general. On the other hand, the basis $V_k$ constructed by TT-sGMRES is non-orthogonal by construction as we perform only an incomplete orthogonalization. Therefore, the truncation in line~\ref{TT-GMRES_orth} only affects the local orthogonality of $V_k$.

In our extensive numerical testing, we experimented with different parameters of the form~\eqref{eq:eta_selection}, possibly including the conditioning of the basis at the denominator as well. However, it turned out that in our context it is good practice to not truncate the vector $\widetilde v_k$ in line~\ref{TT-GMRES_round12} of Algorithm~\ref{alg:sTT-GMRES} (or, equivalently, using a very small $\nu_k$). Indeed, to have a reliable sketching procedure, the update of $W_k$ in line~\ref{TT-sGMRES_updatesketching} should not involve any truncated quantities so that the computation of $y_k$ in~\eqref{eq:pseudo_inv} is coherent with the original, sketched least squares problem~\eqref{eq:sketched_argmin} and not related to a \emph{nearby} problem. See also section~\ref{sec:reorthogonalization} for a similar discussion in case of whitening.

On the other hand, to maintain the TT-ranks of the basis vectors under control, along with selecting small values of $\ell$ (see section~\ref{sec:reorthogonalization}), we perform a 
%rather aggressive
truncation step in line~\ref{TT-GMRES_round22} of Algorithm~\ref{alg:sTT-GMRES}. In particular, the simple strategy of using a constant tolerance $\eta_k\equiv\eta $, for large $\eta$, seems to provide the best trade-off between efficiency (the TT-ranks remain small) and rate of convergence (no remarkable delays have been observed).
% {\color{red} This is probably due to the fact that the low-rank truncations involve only the basis $V_k$ whereas computing $y_k$ as the solution to~\eqref{eq:pseudo_inv} is equivalent to performing an explicit projection with the sketched quantity $S\mathcal AV_k$.}

There are a few cases, in particular when dealing with preconditioned GMRES, that 
we discuss in detail in section~\ref{sec: preconditioning}, where this 
truncation policy is not enough to maintain the TT-rank under control. When 
this happens, we introduce a further parameter \texttt{maxrank} and 
in the truncation phase we use it as a cap on the TT-ranks 
of the basis vectors. This can be done easily within the TT-SVD (performing 
truncated SVDs in all modes) as well as in the randomized schemes that we 
discuss in section~\ref{Randomized approximation of sums in TT-format}.
This action may cause the generated subspace to deviate from the Krylov subspace, losing some theoretical guarantee over the convergence. However, this does not necessarily imply that convergence is lost. For instance, our experiments show that this strategy is very effective when the application of $\mathcal{A}$ leads to an excessive 
growth of the ranks. 
Most importantly, there is no loss of 
accuracy in the projected and true solution, because we ensure that the 
action of the operator is sketched \emph{before} performing the 
rounding. 

% {\color{blue}COMMENTO SU MAXRANK
% The cost of the algorithm highly depends on the ranks of the $v_i$. While the randomized techniques we proposed generally mitigate the growth of these ranks, they may not always suffice. To prevent resource constraints, we recommend setting maximum TT-ranks and truncating tensors once this threshold is surpassed. This action may cause the space we generate to deviate from the Krylov subspace, losing some theoretical guarantee. However, this does not necessarily imply that convergence is lost. For instance, our experiments show that this strategy is very effective when the application of $A$ leads to an unjustified growth of the ranks. 

% }
% { \color{red} Leo: penso che alla fine ci siamo convinti 
% che questa sia sostanzialmente una cattiva idea. Spenderei 
% qualche riga a spiegare perché.}
% \rr{D: non lasciarlo in sottosezione separata.}
% { \color{red}
%   \begin{itemize}
%       \item Riportate il discorso di un opzionale maxrank, 
%         da usare poi nella sezione preconditioning.
%   \end{itemize}
% }

\subsection{Randomized approximation of sums in TT-format}\label{Randomized approximation of sums in TT-format}

As already mentioned, the rounding procedure and the partial orthogonalization
in lines \ref{TT-GMRES_orth}--\ref{TT-GMRES_round22} of Algorithm~\ref{alg:sTT-GMRES}
allow us to mitigate the growth of the TT-ranks due to performing linear combinations of basis vectors. 
The most immediate way 
to implement this operation is to perform 
a rounding after each summation 
in line~\ref{TT-GMRES_orth}. However, this strategy would lead to computing up to $\ell$ extra rounding steps with an excessive increment in the computational efforts. A similar observation applies to the final 
reconstruction of the solution vector in line \ref{alg:final_line22}. 

In this section, we propose to exploit recently developed randomization techniques to reduce 
these costs. Our approach builds upon the algorithms described in \cite{rounding_TT,kressner2023streaming}.
These algorithms are generalizations of randomized low-rank matrix approximation schemes to the tensor realm and provide a significant reduction in computation compared to deterministic algorithms.
In particular, they are particularly effective for rounding or approximating sums of multiple tensors.

The standard deterministic algorithm for TT-rounding is the TT-SVD \cite{TT_decomposition} and requires first to iteratively orthogonalize the TT-cores of the input TT-format. Other approaches incorporating randomization have been proposed, such as the \textit{Randomize-then-Orthogonalize} in \cite{rounding_TT}, which circumvents this orthogonalization step by applying the randomized SVD algorithm~\cite{randomproof} to unfoldings of the full tensor and leveraging the TT-format through the use of Gaussian TT-DRMs (see Definition \ref{defi: TT-DRM}), or a two-sided variant based on generalized Nystr\"om \cite{rounding_TT}. 
The latter 
% An algorithm called \textit{Randomize-then-Orthogonalize} in \cite{rounding_TT} circumvents this orthogonalization step by applying the randomized SVD algorithm~\cite{randomproof} to unfoldings of the full tensor and leveraging the TT-format through the use of Gaussian TT-DRMs (see Definition \ref{defi: TT-DRM}). 
%In~\cite{rounding_TT}, a two-sided variant that uses the generalized Nystr\"om method with two Gaussian TT-DRMs in place of the randomized SVD scheme is also presented.
%The two-sided variant 
has been extended in \cite{kressner2023streaming} to general sketchings, 
and is the algorithm that we will exploit in this work. 
Crucially, the implementation presented in \cite{kressner2023streaming}, called Streaming Tensor-Train approximation (STTA), has the advantage of being \emph{streamable}, namely it requires to operate with the tensor $\mathcal{A}$ only once. This feature will be particularly important in our setting,
as shown later. 

\begin{defi}[Random Gaussian TT-Tensor] \label{defi: TT-DRM} Given a set of target TT-ranks $\{\ell_k\}$, a random Gaussian TT-tensor $\mathcal{L}\in  \mathbb{R}^{n_1\times \dots \times n_d}$ is such that each core tensor $\mathcal{T}_{\mathcal{L},k}\in \mathbb{R}^{\ell_{k-1}\times n_k\times \ell_{k}}$ is filled with random, independent, normally distributed entries with mean 0 and variance $1/(\ell_{k-1} n_k \ell_{k})$ for $1\leq k \leq d$.
\end{defi}

The strength of TT-DRMs is in their ability to reduce the cost of computing partial contractions.  In particular, the $\mu$th right partial contraction of a TT-tensor $\mathcal{T}\in \mathbb{R}^{n_1\times \dots \times n_d}$ of ranks $t_1, \dots, t_{d-1}$ with $\mu$th right interface matrix $C_{>\mu}$ and a Gaussian TT-DRM $\mathcal{R}\in \mathbb{R}^{n_1 \times \dots \times n_d}$ of ranks $r_1, \dots, r_{d-1}$ with $\mu$th right interface matrix $X_{>\mu}$ is the $t_\mu \times r_\mu  $matrix $R_\mu = C_{> \mu}^T X_{>\mu}$. Analogously the $\mu$th left partial contractions of $\mathcal{T}$ and a Gaussian TT-DRM $\mathcal{L}\in \mathbb{R}^{n_1\times \dots\times  n_d}$ and ranks $\ell_1, \dots, \ell_{d-1}$ is the $\ell_{\mu}\times t_\mu$ matrix $L_{\mu} = Y_{\leq \mu}^TC_{\leq \mu} $.

Partial contractions are particularly appealing objects as they can be computed by exploiting the TT structure of the problem, making the computations of the sketchings very cheap. Moreover, having the partial contractions at hand is sufficient to recover the STTA of a tensor. 

The STTA algorithm consists of three phases: the generation phase, the sketching phase, and the recovery phase. In the generation phase, we draw the sketchings, specifically Gaussian TT-DRMs in this case. During the sketching phase, we compute the partial contractions mentioned above. 
Finally, in the recovery phase, we recover the STTA approximant. Below is a summary of the fundamental steps. For more details, please refer to \cite{kressner2023streaming}.

Given a tensor $\mathcal{T}\in \mathbb{R}^{n_1\times \dots \times n_d}$ in TT-format, with ranks $t_1, \dots, t_{d-1}$ and target ranks $r_1, \dots, r_{d-1}$, the STTA algorithm in the generation phase draws random matrices
\[X_{>\mu}\in \mathbb{R}^{(n_{\mu+1}\cdots n_d) \times r_\mu} \quad\text{and}\quad Y_{\leq \mu}\in\mathbb{R}^{(n_1\cdots n_\mu)\times \ell_\mu}, \quad \text{with}\quad \ell_{\mu}>r_{\mu}\hfill, \] then in the sketching phase computes the sketchings 
\[\Psi_\mu = (Y_{\leq \mu-1}^T\otimes I)\mathcal{T}_{\leq \mu} X_{>\mu} \quad\text{and}\quad \Omega_{\mu} = Y_{\leq \mu}^T \mathcal{T}_{\leq \mu} X_{>\mu},\]
and finally forms the right unfoldings of the TT-cores $\widehat{C}_\mu$ as
\[\widehat{C}_\mu^{R} = \Omega_{\mu-1}^{\dagger}\Psi_{\mu}.\]
A possible way to construct the sketching matrices $X_{>\mu}$ and $Y_{\leq\mu}$ is to use respectively the right and left interface matrices of two Gaussian TT-DRMs of appropriate size.

These steps describe how to compute the STTA approximation of a tensor. To compute the STTA approximant of a linear combination of tensors $a_1\mathcal{T}^{(1)}+\dots + a_s\mathcal{T}^{(s)}$, first use the same DRMs to sketch each $\mathcal{T}^{(i)}$ obtaining the $\Psi_{\mu}^{(i)}$ and $\Omega_{\mu}^{(i)}$. Next, compute the linear combinations $\Psi_\mu =  a_1\Psi_\mu^{(1)}+\dots + a_s\Psi_\mu^{(s)}$ and $\Omega_\mu =  a_1\Omega_\mu^{(1)}+\dots + a_s\Omega_\mu^{(s)}$. Finally, proceed as described above to recover the final approximant.

 The STTA algorithm can be exploited in TT-sGMRES during the orthogonalization phase, to compute the weighted sum in line~\ref{TT-GMRES_orth} of Algorithm~\ref{alg:sTT-GMRES}, and in line~\ref{alg:final_line22} to compute the final solution.
In particular, since we only need the sketched matrices $\Omega_\mu^{(v_k)}$ and $\Psi_\mu^{(v_k)}$ of each basis vector $v_k$ to form the final solution $x_k$ using STTA, we can get rid of the basis vectors that are no longer needed in the incomplete orthogonalization and store only their sketches. This is particularly beneficial in situations where memory constraints pose a challenge. This means that we can exploit the full potential of the incomplete orthogonalization also in terms of storage demand while avoiding the possible extra costs coming from a two-pass strategy.

In practice, we have implemented the rounding schemes proposed in 
\cite{kressner2023streaming}, and obtained two routines, called 
\textsc{STTA\_sketch} and \textsc{STTA\_recover}, that perform the 
following actions:
\begin{description}
    \item[\normalfont\textsc{STTA\_sketch}] takes as input 
      a tensor $\mathcal T$ and $X, Y$ as described above, 
      and computes the corresponding sketches 
      $\Psi_\mu$ and $\Omega_\mu$.
    \item[\normalfont\textsc{STTA\_recover}] takes as input the sketches 
      $\Psi_\mu$ and $\Omega_\mu$ (resp. a linear combination 
      of sketchings) and reconstruct an approximation to 
      the original tensor 
      $\mathcal T$ (resp. the linear combination of the tensors). 
\end{description}
Throughout the algorithm, we assume that the tensors $X, Y$ have been 
chosen at the beginning, with suitable dimensions $r_\mu, \ell_\mu$, 
which we discuss in further detail in section~\ref{sec:alltogether}.
We are not able to recommend a choice for these 
parameters that is suitable for all 
cases; in the algorithms we let the user provide the values of these 
parameters.

\subsection{Incomplete orthogonalization, restarting, and whitening}
\label{sec:reorthogonalization}

From a computational point of view, being able to perform only a local orthogonalization in line~\ref{TT-GMRES_orth} of Algorithm~\ref{alg:sTT-GMRES} is key to attain a competitive solver. However, choosing a suitable value of $\ell$, the scalar that controls the number of vectors to orthogonalize the newly computed basis vector against, is not straightforward. This is a common issue also in the case of truncated Krylov methods for standard linear systems of equations; see, e.g.,~\cite{Saad_truncatedGMRES}.

In our context, employing smaller values of $\ell$ not only decreases the cost of the orthogonalization step itself, thanks to fewer orthogonalizations to perform, but it also induces smaller TT-ranks in the result by reducing the number of tensor sums. 
This means that adopting a very small $\ell$ has an impact on the whole solution procedure and is extremely beneficial in reducing the computational efforts devoted to every operation
involving the basis vectors in TT-format.
In most of our experiments we select $\ell=1$ obtaining a very successful solution process; see section~\ref{sec: numerical experiments}.

If selecting a small $\ell$ looks very appealing from a computational point of view, such a selection most likely leads to a basis $V_k$ which is terribly ill-conditioned. In~\cite[Section 5.3]{nakatsukasa2021fast} the authors suggest to restart the iterative scheme whenever a too ill-conditioned basis $V_k$ is detected. In particular, if at iteration $m$, $V_m$ turns out to be (close to) singular, we may construct 
the residual vector $r_m=b-\mathcal{A}x_m$, and restart the TT-sGMRES iteration using $r_m$ as new initial residual vector in line~\ref{sTTGMRES_line1} of Algorithm~\ref{alg:sTT-GMRES}. Even though this machinery may help in reducing the impact of working with an ill-conditioned basis, it can potentially lead to important delays in the convergence of the overall solution process. In practice, we have never needed to employ this strategy 
in our numerical experiments. Moreover, in~\cite{GuettelSimunec2023} it has been observed that having an ill-conditioned $V_k$ is not the primary cause of the possible numerical instabilities of sGMRES. Therefore, we do not adopt any restarting strategy in our numerical examples.

A different approach to stabilize sketched Krylov methods is the so-called \emph{whitening}, namely performing an explicit full orthogonalization of the sketched basis $SV_k$.
This inexpensive procedure has a rather important impact in our context as it allows us to rewrite the minimization problem~\eqref{eq:pseudo_inv} in a different way, reminiscent of the projected formulation~\eqref{eq:minres2} of (standard) GMRES.
In particular, in~\cite{PalittaSchweitzerSimoncini2023} a sketched Arnoldi relation has been derived in the context of Krylov approximations to matrix function evaluations. Let $V_k$ be constructed by a truncated Arnoldi scheme for which the Arnoldi relation~\eqref{eq:ArnoldiRel} holds true. Moreover, let $Q_kT_k=SV_k$ be the skinny QR factorization of the sketched basis $SV_k$ and
$$SV_{k+1}=[Q_k,q_{k+1}]\begin{bmatrix}
T_k & t_{k+1} \\
0 & \tau_{k+1}\\
\end{bmatrix}.
$$
Then, we can write 
\begin{equation}\label{eq:sketchedArnoldi}
    SA\widehat V_k=S\widehat V_k(\widehat H_k+\widehat he_k^T)+h_{k+1,k}S\widehat v_{k+1}e_k^T=S\widehat V_{k+1}\begin{bmatrix}
\widehat H_k+\widehat he_k^T\\
[0,\ldots,0,h_{k+1,k}]\\
    \end{bmatrix},
\end{equation}
where $\widehat V_{k+1}=[\widehat v_{1},\ldots,\widehat v_{k+1}]=V_{k+1}T_{k+1}^{-1}$, $\widehat H_k=T_kH_kT_k^{-1}$, and $\widehat h=t_{k+1}h_{k+1,k}/\tau_k$; see~\cite[Equation 9]{PalittaSchweitzerSimoncini2023}. 
Even though the transformed basis $\widehat V_{k+1}$ is not explicitly available, it is important to notice that this is orthogonal with respect to the sketched inner product $S^TS$, namely $\widehat V_{k+1}^TS^TS\widehat V_{k+1}=I$.

Thanks to~\eqref{eq:sketchedArnoldi} and the $S^TS$-orthogonality of $\widehat V_k$, the minimization problem~\eqref{eq:sketched_argmin} can be reformulated as
\begin{align}\label{eq:sketched_argmin2}
    y_k = &\;\mathrm{argmin}_y \|SAV_k y - Sr_0\|_2=\mathrm{argmin}_y \|SAV_kT_{k}^{-1}T_k y - Sr_0\|_2\notag\\
    =&\;\mathrm{argmin}_{y=T_k^{-1}z} \|S\widehat V_{k+1}z - Sr_0\|_2\notag\\
    =&\;
    \mathrm{argmin}_{y=T_k^{-1}z} \left\|\begin{bmatrix}
\widehat H_k+\widehat he_k^T\\
[0,\ldots,0,h_{k+1,k}]\\
    \end{bmatrix}z - \beta e_1\right\|_2,\quad \beta=\|Sr_0\|_2.
\end{align}
If the vector $y_k$ is computed as above, the sketched norm of the residual vector associated to the solution $x_k=x_0 +V_k y_k$, namely $r_k=b-Ax_k$, can be cheaply computed as 
\begin{equation}\label{eq:sketched_residualnorm}
\|r_k\|=\|S(A V_k y_k-r_0)\|=\|S(A \widehat V_k z_k-r_0)\|=\left\|\begin{bmatrix}
\widehat H_k+\widehat he_k^T\\
[0,\ldots,0,h_{k+1,k}]\\
    \end{bmatrix}z_k - \beta e_1\right\|_2.    
\end{equation}

We would like to mention that, to the best of our knowledge, the derivations above result to be new, even though they come from a straightforward combination of the original sGMRES scheme from~\cite{nakatsukasa2021fast} and the sketched Arnoldi relation presented in~\cite{PalittaSchweitzerSimoncini2023}.

If one wanted to adopt whitening, the only operations to change in Algorithm~\ref{alg:sTT-GMRES} would be the computation of $y_k$ in line~\ref{sTT-GMRES_compute_y} and the residual norm evaluation in line~\ref{sTT-GMRES_compute_resnorm}. Moreover, the storage of the matrix $W_k$ would be no longer necessary whereas the updating of the skinny QR factorization of $SV_k$ would have to be introduced. 

Even though it has been shown that whitening is an extremely beneficial practice in contexts like matrix function approximations~\cite{PalittaSchweitzerSimoncini2023} and the numerical solution of matrix equations~\cite{PalittaSchweitzerSimoncini2023bis}, we must mention that it does present some peculiar drawbacks in our framework. In particular, the computation of the coefficients collected in the matrix $H_d$ takes place before truncating the current basis vector $\widetilde v$ in line~\ref{TT-GMRES_round22} of Algorithm~\ref{alg:sTT-GMRES}. On the other hand, the sketching $S$ is applied to $v_{k+1}$, the truncated (and normalized) version of $\widetilde v$. $Sv_{k+1}$ is then used to update the skinny QR of $SV_{k+1}$ and thus obtain the coefficients in $T_{k+1}$ necessary for computing the quantities involved in the projected problem~\eqref{eq:sketched_argmin2}. 
As it turned out from our vast numerical testing, this discrepancy in the construction of $H_k$ and $T_k$ may lead to a disagreement between the actual sketched residual norm $\|SAV_ky_k-Sr_0\|$ and its computed value on the right-hand side of~\eqref{eq:sketched_residualnorm}, whenever $y_k$ is computed as in~\eqref{eq:sketched_argmin2}. We did not observe such a trend when computing $y_k$ by~\eqref{eq:pseudo_inv}. Indeed, the use of the pseudoinverse of $SAV_k$ is equivalent to performing an explicit projection without relying on the sketched Arnoldi relation~\eqref{eq:sketchedArnoldi}. Therefore, in all the experiments reported in section~\ref{sec: numerical experiments} the vector $y_k$ is computed by~\eqref{eq:pseudo_inv}.

\subsection{Building the final solution }
\label{sec:final-solution}

The final step of the TT-sGMRES algorithm is the computation of the solution $x_k=x_0+ V_ky_k = x_0+\sum_i^k v_i [y_k]_i$. For this task, we propose to use the STTA algorithm.

Compared with the classic way to perform this linear combination (adding one term at a time and rounding after each addition), this algorithm offers several advantages, some of which we have already described at the beginning 
of section~\ref{Sketched TT-GMRES}. In particular, this strategy has a lower computational costs and avoids the storage of the basis.
Another advantage is that when the basis $V$ is not orthogonal, 
possibly badly conditioned, the classic procedure may face numerical cancellation. On the other hand, our results show that STTA is not affected by this undesirable issue. There is, however, a drawback in using STTA. Indeed, this strategy requires knowing in advance the numerical TT-rank of the solution, or at least an overestimate thereof, which is not available in general. For the moment, we lack valid automatic strategies for estimating the TT-rank of the final solution and in our routines we rely on a user-provided value. That said, for many problems of interest, the TT-ranks of the solution are very low, even lower than those of a single $v_i$, so that any reasonable heuristic could work.

\subsection{Putting it all together}
\label{sec:alltogether}

In Algorithm~\ref{alg:enhanced sTT-GMRES} we report the TT-sGMRES pseudocode enhanced with all the tools and considerations discussed in the previous sections. In particular, as mentioned in section~\ref{sec:truncation}, we refrain from performing any low-rank truncation after the application of $\mathcal{A}$ in line~\ref{application_of_A_TT-sGMRES} whereas we employ a rather large, constant value $\eta$ ($\eta$ is 
either $0.1$ or $0.3$ in our experiments, and we choose it to ensure 
that the prescribed tolerance is 
reached) in the truncations in line~\ref{rounding_GS_enhancedTTsgmres}. Moreover, any linear combinations involving the basis TT-vectors (line~\ref{STTArecover1} and~\ref{alg:final_line222}) is carried out by the \textsc{STTA\_recover} routine described in section~\ref{Randomized approximation of sums in TT-format}. To this end, we compute the sketch of the newly defined basis vector $v_{k+1}$ by \textsc{STTA\_sketch} in line~\ref{STTArecover2}. The parameter $\ell_\mu$ for the STTA algorithm (the oversampling) is set to $20$. 

The number of rows of the sketch $S$ for the TT-sGMRES method is based on the maximum number of iterations. If the user specifies a maximum number \texttt{maxit}, the number of rows of $S$
is chosen as twice that number.
Optionally, in our code we allow to further tweak this parameters, 
or to specify a custom sketching $S$. 

\begin{remark}
    In the pseudocode of Algorithm~\ref{alg:enhanced sTT-GMRES}
    we use the routines
    \textsc{STTA\_Sketch} and \textsc{STTA\_Recover} to perform 
    the partial reorthogonalization. This is useful especially for sizable values of $\ell$. However, in our experiments we often choose $\ell = 1$, for which it 
    is instead preferable to maintain in memory the last vector 
    and perform the reorthogonalization and rounding explicitly in the 
    TT-format. In our implementation we let the user choose between the 
    two stratgies.
\end{remark}

%\begin{algorithm} 
    %\begin{algorithmic}[1]
    %\Procedure{Enhanced sTTGMRES}{$A,b,x_0,S, k, \tau$}
        %\State $\eta \gets 10^{-1}$ \Comment{Adjustment for the %tolerance, see \ref{sec:truncation}}
     %   \State $r_0 \gets b - A x_0$
     %   \State $v_1 \gets r_0 / \| r_0 \|_2$
    %    \State $r_{0,S} \gets S r_0$ 
   %        \Comment{Sketch of the initial residual}
  %      \State $\mathrm{res} \gets \norm{r_{0,S}}$
 %       \While{$\mathrm{res} > \tau $}
%            \State $w_i \gets \Call{Round}{A v_i, \eta \tau}$ %\label{alg:Av}
  %          \State $w_{i,S} \gets Sw_i$
  %          \State $w_i \gets w_i - 
  %             \sum_{j = 1}^{\min(k, i)}  \langle 
  %               v_{i-j+1}^T, w_i \rangle_S\ 
  %               v_{i-j+1}$
  %            \Comment{see \ref{sec:reorthogonalization}}
  %          \State $v_{i+1} \gets 
  %            \Call{Round}{w_i / \norm{w_i}_2, \eta \tau}$
  %            \label{alg:rounding}
  %          \State $W \gets \begin{bmatrix}
  %              w_{1,S} & \dots & w_{i,S}
  %          \end{bmatrix}$
  %          \State $y \gets W^\dagger r_S$
  %          \State $\mathrm{res} \gets \| r_{0,S} - W y \|$
 %       \EndWhile
 %       \State \Return% 
 %         \Call{Round}{$y_1 v_1 + \ldots + y_i v_i, \tau$} %\label{alg:reconstruction}
          %\Comment{Sum computed by sketching, see \ref{sec:final-%solution}}
%        \EndProcedure
%    \end{algorithmic}
%    \caption{Sketched TT-GMRES algorithm for the 
%  linear system $Ax = b$, with tolerance $\tau$, 
%  starting guess $x_0$, 
%  sketching $S$, and orthogonalization parameter $k$.}
%    \label{alg:enhanced sTT-GMRES}
%\end{algorithm}

\begin{algorithm}[t]%\caption{TT-GMRES}
\begin{algorithmic}[1]
%\setstretch{1.2}
\smallskip
\Statex \textbf{Input:} Tensor $\mathcal A\in\mathbb{R}^{n_1\times\ldots\times n_d}$, right-hand side $b$, initial guess $x_0$ in TT-format, maximum basis dimension $\texttt{maxit}$, tolerance $\texttt{tol}$, sketching $S$, incomplete orthogonalization parameter $\ell$, rounding threshold $\eta$.  
\Statex \textbf{Output:} Approximate solution $x_k$  such that $\|S(\mathcal Ax_k-b)\|\leq \|Sb\| \cdot \texttt{tol}$
\smallskip

\State Set $r_0=b-\mathcal Ax_0$, $\beta=\|r_0\|$ $V_1=v_1=r_0/\beta$, $\beta^{[S]} = \| S b \|$ \label{sTTGMRES_line1}, $W_0 = []$%, $Q_1=SV_1/\|SV_1\|$, $T_1=\|SV_1\|$
\State $[\Phi^{(1)}, \Psi^{(1)}] = 
  \Call{STTA\_Sketch}{v_1, X, Y}$,
 \For{$k=1,\ldots,\texttt{maxit}$}
\State Compute $\widetilde v=\mathcal A v_k$ \label{ApplyATT-sGMRES}% \Comment{Adjustment for the tolerance, see \ref{sec:truncation}}
% \State Update $S\mathcal{A}V_k=[S\mathcal{A}V_{k-1},S\widetilde v]$\label{application_of_A_TT-sGMRES}
%\State Compute $\widetilde v^{[S]} = $
\State Update $W_k=[W_{k-1},S \widetilde v]$\label{application_of_A_TT-sGMRES}
\For{$i=\max\{{1,k-\ell+1,\ldots,k}\}$}
\State 
  %Set $\widetilde v=\widetilde v- v_ih_{i,k}$, where $h_{i,k}= \widetilde v^Tv_i$ 
  Set $h_{i,k}= \widetilde v^T v_i$
  \Comment{Only $\ell$ previous vectors are kept in memory}
 % {\color{red} DP: usa STTA\_recover}
 \label{STTArecover1}%              \Comment{see \ref{sec:reorthogonalization}}
\EndFor
\For{$\mu = 1, \ldots, d$}
\State Set $\widetilde \Phi_\mu = h_{1,k} \Phi_\mu^{(1)} + \ldots  + h_{1,k} \Phi_\mu^{(k)}$ and $\widetilde \Psi_\mu = h_{1,k} \Psi_\mu^{(1)} + \ldots  + h_{1,k} \Psi_\mu^{(k)}$
\EndFor
\State Set $\widetilde v=\Call{STTA\_Recover}{\widetilde \Phi, \widetilde \Psi, \eta\cdot \texttt{tol}}$\label{rounding_GS_enhancedTTsgmres} %\Comment{Adjustment for the tolerance, see \ref{sec:truncation}}
\State Set $h_{k+1,k}=\|\widetilde v\|$ and $v_{k+1}=\widetilde v/h_{k+1,k}$
% \State {\color{red} DP: skectha il nuovo vettere $v_{k+1}$ con STTA\_sketch} 
\State Compute  $[\Phi^{(k+1)}, \Psi^{(k+1)}] = \Call{STTA\_Sketch}{v_{k+1}, X, Y}$
\label{STTArecover2}
%\State Update skinny QR: $Q_{k+1}T_{k+1}=[SV_{k},Sv_{k+1}]$ \label{white:line12}
%\State Update $\widehat H_k=T_kH_kT_k^{-1}$ and set $\widehat h=(h_{k+1,k}/\tau_k)t_{k+1}$ \label{white:line22}
\State Compute $y_k$ as the solution to~\eqref{eq:pseudo_inv}%\label{sTT-GMRES_compute_y2}
%\State Compute $\|r_k\|$ as in~\eqref{eq:sketched_residualnorm}
\If{$\|W_ky_k-Sr_0\|\leq \beta^{[S]}\cdot\texttt{tol}$}%\label{sTT-GMRES_compute_resnorm2} 
\State Go to line~\ref{alg:final_line222}
\EndIf
\State Set $V_{k+1}=[V_k,v_{k+1}]$ 
\EndFor
\For{$\mu = 1, \ldots, d$}
\State Set $\widetilde \Phi_\mu = [y_{k}]_1 \Phi_\mu^{(1)} + \ldots  + [y_{k}]_k \Phi_\mu^{(k)}$ and $\widetilde \Psi_\mu = [y_{k}]_1 \Psi_\mu^{(1)} + \ldots  + [y_{k}]_k \Psi_\mu^{(k)}$
\EndFor
\State Set $x_k = \Call{STTA\_Recover}{\widetilde \Psi, \widetilde \Phi, \mathrm{tol}}$
%\State Set $x_k=\Call{Round}{x_0+V_ky_k,\mathrm{tol}}$  {\color{red} DP: usa STTA\_recover}
\label{alg:final_line222}%           \Comment{Sum computed by sketching, see \ref{sec:final-solution}}

\end{algorithmic}    \caption{Sketched TT-GMRES (TT-sGMRES)}
    \label{alg:enhanced sTT-GMRES}

\end{algorithm}

\section{Preconditioning}\label{sec: preconditioning} 

It is well-known that, to get a fast rate of convergence in terms of number of iterations, Krylov methods require preconditioning in general. This applies to our TT-sGMRES scheme as well. However, due to the peculiarity of our framework, preconditioners for~\eqref{eq:main} may pose
further challenges with respect to preconditioninig operators for standard linear systems. Indeed, in addition to be effective in reducing the number
of iterations at a reasonable computational cost, the preconditioner operator must not
dramatically increase the rank of the current basis vector. Otherwise, the cost of all the remaining operations in TT-sGMRES would increase possibly jeopardizing the gains coming from running fewer iterations. A similar scenario holds 
for standard TT-GMRES as well. 

Note that, in principle, thanks to the incomplete orthogonalization we perform, TT-sGMRES 
is less penalized than the standard TT-GMRES~\cite{TT-GMRES}
if a large number of iterations to converge is needed. 
Nevertheless, for several practical problems (for instance the 
ones arising from PDEs, where the condition number of the problem 
grows with the problem dimension), preconditioning is essential to ensure 
convergence in a reasonable amount of time. 

Few options for preconditioning tensor equations of the form~\eqref{eq:main} are available in the literature. In~ \cite{Precond1}, a low-rank approximation to $\mathcal{A}^{-1}$ is employed as preconditioner for~\eqref{eq:main}. Exponential sums have been 
proposed in \cite{hackbusch2015hierarchical,hackbusch2019computation, TTGMRES2,rohrig2023performance}. 

The main limitation when dealing with preconditioning
in tensor Krylov methods is that the operator 
$\mathcal A \mathcal P^{-1}$ is usually of a much higher 
tensor rank than $\mathcal A$, and therefore induces
a much faster rank growth in the basis. Hence, 
even if the number of iterations necessary for 
convergence can be greatly reduced, this does 
not necessarily correspond to a reduction in 
computational cost. In the next section, we discuss 
how sketching can be helpful in this context as well, 
by limiting the maximum TT-rank that can be 
reached in the GMRES basis. 

We could consider left or right preconditioning, or both at once. We choose to 
only discuss right preconditioning because it ensures that the residual of the 
preconditioned problem and of the original one coincide. In a nutshell, assuming 
the availability of a preconditioner $\mathcal P$, right preconditioning modifies 
lines~\ref{ApplyATT-sGMRES} and \ref{alg:final_line222} in Algorithm~\ref{alg:enhanced sTT-GMRES} as follows:
\begin{align*}
    \widetilde v & =  \mathcal A \mathcal P^{-1} v_k, \\ 
    x_k & = \mathcal P^{-1} \left[ \textsc{STTA\_Recover}(\widetilde \Psi, \widetilde \Phi, \mathrm{tol}) \right].
\end{align*}
As we discuss in section~\ref{sec:maxrank}, 
this does not always lead to better performances even when 
the preconditioner works nicely, and some extra care needs to be taken to avoid an excessive rank growth. In particular, it turned out that coupling preconditioning 
with a ``maximum rank'' rounding step and sketching often leads to competitive results. 

\subsection{Exponential sum preconditioning}

 In this work, we have considered preconditioners based on \emph{exponential sums}, that are often suitable for problems arising from PDEs; see, e.g., \cite{hackbusch2015hierarchical,hackbusch2019computation, TTGMRES2,rohrig2023performance}.
 In order to construct such preconditioner, it is first necessary to split the operator 
$\mathcal A$ into the form 
\[
  \mathcal A = \widehat{\mathcal A} + \bigoplus_{i = 1}^d A_i, 
\]
where the second term is the dominant part of the operator and $\bigoplus$ denotes a Kronecker sum, that is
\[
  \bigoplus_{i = 1}^d A_i = 
  A_d \otimes I \otimes \ldots \otimes I + \ldots + 
  I \otimes \ldots \otimes I \otimes A_1. 
\]
The above is a summation of $d$ terms, each with a single 
entry in the Kronecker product different from the identity, 
which form a commutative family. 
% {\color{red}DP: do we really need $d$ terms in the preconditioner? $d$ is the number of mode of $\mathcal{A}$.} 
Then, we precondition by considering $\mathcal P$ such 
that $\mathcal P^{-1} \approx (\bigoplus_{i = 1}^d A_i)^{-1}$. Instead of computing explicitly such $\mathcal P$, we directly write $\mathcal P^{-1}$. 
 To accomplish this, we rely on
 exponential
 sums, that is we determine an approximant for the inverse 
 function $\frac{1}{z}$ of the form 
 \[
   \frac{1}{z} = \sum_{j = 1}^{\zeta} \alpha_j e^{-\beta_j z} =: E_{\zeta}(z), 
 \]
 where $\zeta$ is a positive integer, and such that the approximation 
 is accurate over the spectrum (or better, over the field of values) of $\bigoplus_{i = 1}^d A_i$. Then, we 
 consider 
 \[
  \mathcal P^{-1} := E_\zeta\left(\bigoplus_{i=1}^d A_i\right) = 
   \sum_{i = 1}^\zeta \alpha_i 
     \bigotimes_{j = 1}^d e^{-\beta_i A_j}. 
 \]
 % {\color{red}DP: double check $i$ and $j$ above.}
 
 In particular, applying $\mathcal{P}^{-1}$ to a 
 tensor $\mathcal X$ requires to sum $\zeta$ tensors, obtained by performing 
 $j$-mode multiplications with $e^{-\beta_i A_j}$ for all $j$. Since 
 $j$-mode multiplications do not increase the TT-rank, 
 applying this preconditioner generally increases
 the TT-ranks 
 of $\mathcal X$ by a factor of (at most) $\zeta$.

 The difficulty in designing a preconditioner in this 
 class lies in determining the coefficients 
 $\alpha_i, \beta_i$. In this work we rely on the 
 procedure described in \cite{TTGMRES2}; we refer 
 the interested reader to  \cite{hackbusch2019computation} and \cite[Appendix D]{hackbusch2015hierarchical} 
 for an in-depth overview.
 Determining the optimal $\alpha_i$, $\beta_i$ is often 
 challenging even when the spectrum is real and 
 known a-priori (see \cite{hackbusch2019computation}); hence, 
 we often prefer to rely on suboptimal approximations 
 recovered from integral representations of $\frac{1}{z}$ (as done in \cite{TTGMRES2}). 
 It is worth noting that another approach to preconditioning this class of problems involves techniques based on tensor Sylvester equations, such as those presented in \cite{casulli2023tensorized}.

\subsection{Sketching and bounded rank roundings}
\label{sec:maxrank}

We note that several techniques discussed in the previous sections (e.g., incomplete reorthogonalization) might become less relevant when using 
a good preconditioner as this leads to convergence 
in a small number of steps, in general. On the other hand, preconditioning often leads to fast rank growth,
possibly making the overall solution process impractical. To mitigate this annoying side-effect, we propose to rely on 
a low-rank rounding step 
of the basis with a prescribed maximum rank. This gives little control over the truncation accuracy, making the analysis of the method even trickier. In particular, the distance between truncated and original (not truncated) quantities cannot be quantified in general. %{\color{red}DP: ho tolto la frase su Krylov che secondo me apriva a troppe domande}
However, sketching-based GMRES still works fine in practice and the maximum-rank rounding often leads to important computational advantages. Nevertheless, we must mention that this rounding may induce a slightly larger (but faster) number of iterations when compared to the scenario where this is not performed.

To implement the maximum-rank rounding, when we call the rounding procedure in line \ref{rounding_GS_enhancedTTsgmres}, we enforce that the TT-rank of $v_{k+1}$
cannot be larger than a maximum prescribed value $r_{\max}$ (component-wise). The choice 
of this $r_{\max}$ is arbitrary and the optimal value problem dependent: smaller ranks 
correspond to faster iterations but slower convergence, whether higher ranks lead to 
fewer iterations but with a higher computational cost per iteration. 

The preconditioned variant of TT-sGMRES, that we call TT-sPGMRES, is reported in Algorithm~\ref{alg:Preconditioned_sTT-GMRES}.

\begin{algorithm}[t]

\begin{algorithmic}[1]

\smallskip
\Statex \textbf{Input:} Tensor $\mathcal A\in\mathbb{R}^{n_1\times\ldots\times n_d}$, preconditioning operator $\mathcal{P}$, right-hand side $b$, initial guess $x_0$ in TT-format, $\mathrm{maxit}$ maximum basis dimension, tolerance $\mathrm{tol}$, sketching $S$, incomplete orthogonalization parameter $\ell$, rounding threshold $\eta$, $r_{\max}>0$ maximum rank allowed in the basis vectors.  
\Statex \textbf{Output:} Approximate solution $x_k$  such that $\|\mathcal S(Ax_k-b)\|\leq \mathrm{tol} \cdot \| S b \|$
\smallskip

\State Set $r_0=b-\mathcal A \mathcal P^{-1} x_0$, $\beta=\|r_0\|$ $V_1=v_1=r_0/\beta$, $\beta^{[S]} = \| S b \|$, $W_0 = []$
\State Compute $[\Phi^{(1)}, \Psi^{(1)}] = \Call{STTA\_Sketch}{v_1, X, Y}$
 \For{$k=1,\ldots,\text{maxit}$}
\State Compute $\widetilde v=\mathcal A \mathcal P^{-1}v_k$ \label{precondApplyATT-sGMRES}% \Comment{Adjustment for the tolerance, see \ref{sec:truncation}}
% \State Update $S\mathcal{A}V_k=[S\mathcal{A}V_{k-1},S\widetilde v]$\label{application_of_A_TT-sGMRES}
\State Update $W_k=[W_{k-1},S \widetilde v]$\label{precondapplication_of_A_TT-sGMRES}
\For{$i=\max\{{1,k-\ell+1,\ldots,k}\}$}
\State 
  %Set $\widetilde v=\widetilde v- v_ih_{i,k}$, where $h_{i,k}= \widetilde v^Tv_i$ 
  Set $h_{i,k}= \widetilde v^T v_i$
 % {\color{red} DP: usa STTA\_recover}
 \label{precondSTTArecover1}%              \Comment{see \ref{sec:reorthogonalization}}
\EndFor
\For{$\mu = 1, \ldots, d$}
\State Compute $\widetilde \Phi_\mu = h_{1,k} \Phi_\mu^{(1)} + \ldots  + h_{1,k} \Phi_\mu^{(k)}$ and $\widetilde \Psi_\mu = h_{1,k} \Psi_\mu^{(1)} + \ldots  + h_{1,k} \Psi_\mu^{(k)}$
\EndFor
\State Set $\widetilde v=\Call{STTA\_Recover}{\widetilde \Phi, \widetilde \Psi, \eta\cdot \mathrm{tol}, r_{\max}}$\label{precondrounding_GS_enhancedTTsgmres} %\Comment{Adjustment for the tolerance, see \ref{sec:truncation}}
\State Set $h_{k+1,k}=\|\widetilde v\|$ and $v_{k+1}=\widetilde v/h_{k+1,k}$
% \State {\color{red} DP: skectha il nuovo vettere $v_{k+1}$ con STTA\_sketch} 
\State Compute  $[\Phi^{(k+1)}, \Psi^{(k+1)}] = \Call{STTA\_Sketch}{v_{k+1}, X, Y}$
\label{precondSTTArecover2}
%\State Update skinny QR: $Q_{k+1}T_{k+1}=[SV_{k},Sv_{k+1}]$ \label{white:line12}
%\State Update $\widehat H_k=T_kH_kT_k^{-1}$ and set $\widehat h=(h_{k+1,k}/\tau_k)t_{k+1}$ \label{white:line22}
\State Compute $y_k$ as the solution to~\eqref{eq:pseudo_inv}%\label{sTT-GMRES_compute_y2}
%\State Compute $\|r_k\|$ as in~\eqref{eq:sketched_residualnorm}
\If{$\|W_ky_k-Sr_0\|\leq \beta^{[S]}\cdot\mathrm{tol}$}%\label{sTT-GMRES_compute_resnorm2}
\State Go to line~\ref{alg:precondfinal_line222}
\EndIf
\State Set $V_{k+1}=[V_k,v_{k+1}]$ 
\EndFor
\label{alg:precondfinal_line222}
\For{$\mu = 1, \ldots, d$}
\State Set $\widetilde \Phi_\mu = [y_{k}]_1 \Phi_\mu^{(1)} + \ldots  + [y_{k}]_k \Phi_\mu^{(k)}$ and $\widetilde \Psi_\mu = [y_{k}] \Psi_\mu^{(1)} + \ldots  + [y_{k}]_k \Psi_\mu^{(k)}$
\EndFor
\State Set $x_k = \Call{STTA\_Recover}{\widetilde \Psi, \widetilde \Phi, \mathrm{tol}}$
\State Update $x_k = \mathcal P^{-1} x_k$
%\State Set $x_k=\Call{Round}{x_0+V_ky_k,\mathrm{tol}}$  {\color{red} DP: usa STTA\_recover}
%           \Comment{Sum computed by sketching, see \ref{sec:final-solution}}
\end{algorithmic}    \caption{Preconditioned sketched TT-GMRES (TT-sPGMRES)}
   \label{alg:Preconditioned_sTT-GMRES}

\end{algorithm}

% \subsection{Advantages over AMEn}
% \rr{
% \begin{itemize}
%     \item We can handle matrix-free problems, for instance 
%     exploiting sparsity in the operator, and similar 
%     situations;
%     \item We can apply preconditioners.
% \end{itemize}}

\section{Numerical illustration}\label{sec: numerical experiments}

In this section, we analyze the proposed enhanced TT-sGMRES algorithm through two distinct applications: one involving convection-diffusion PDEs and another arising from Markov chains in performability and reliability 
analysis. We compare its performance against other solvers in 
the TT-format, including TT-GMRES, the vanilla version of TT-sGMRES, 
and AMEn.

A key aspect of the enhanced TT-sGMRES algorithm is that it provides access only to the sketched residual \eqref{eq:sketched_argmin}, which is typically slightly smaller than the actual residual. To ensure fair comparisons, we set the tolerance for TT-sGMRES lower than that of TT-GMRES, guaranteeing that the desired accuracy is consistently achieved across all tested scenarios.

The section is divided into two main blocks, in which we analyze respectively the behaviours of the algorithms without and with preconditioning. Before 
presenting these two blocks experiments in sections~\ref{sec:numexp-unpreconditioned} and \ref{sec:numexp-preconditioned}, respectively, 
we briefly describe the two case studies. 
In all unpreconditioned experiments, the maximum number of iteration
for TT-sGMRES
is set to $200$ (and thus the sketch $S$ has $400$ rows), 
whereas in the preconditioned examples this number is set to $20$ (and $S$ has $40$ rows). 

The code to replicate the numerical experiments in this section 
can be downloaded from \url{https://github.com/numpi/tt-sgmres}. 
It requires MATLAB and the TT-Toolbox \cite{TT-Toolbox}. 

\subsection{Case studies}
\label{sec:case-studies}

Throughout the numerical experiments, we will 
consider two classes of linear systems, that are 
briefly described here. The first arises
from the discretization of a PDE, whereas the 
second stems from the analysis of a high-dimensional 
Markov chain. \medskip 

\subsubsection{A convection-diffusion problem}
\label{sec:case-study-convection-diffusion}
We consider the computation of the 
steady-state for a 
convection-diffusion equation on a $d$-dimensional 
box
\[
  K \Delta u + \langle w ,\nabla u\rangle + f = 0, 
  \qquad 
  u: [-1, 1]^d \to \mathbb R, 
\]
with zero Dirichlet boundary conditions. We
choose the parameters $K = 10^{-2}$ and 
  $w = 10^{-2}\cdot[1,\ldots, 1]\in\mathbb{R}^d$. The source term is chosen 
as $f(x) = e^{-10 \norm{x}^2_2}$. When discretized 
with finite differences this yields the linear
system
\[
  \left( \bigoplus_{i = 1}^d 
    [ L + D_i ] \right) x +  f = 0, 
\]
where $f$ contains the samplings of the source 
term at the grid points, and the matrices 
$L$ and $D_i$ discretize the diffusion 
and convection operators, and 
are defined as follows:
\[
    L = \frac{K}{h^2} \begin{bmatrix}
        -2 & 1 \\
        1 & \ddots & \ddots \\ 
        & \ddots & \ddots & 1 \\
        & & 1 & -2 \\
    \end{bmatrix}, \qquad 
    D_i = \frac{w_i}{h} \begin{bmatrix}
        -1 & 1 \\
        & \ddots & \ddots \\ 
        &  & \ddots & 1 \\
        & &  & -1 \\
    \end{bmatrix}. 
\]
The choice of the source term $f(x, y) = e^{-10 (x^2 + y^2)}$
guarantees that, when represented in tensor form, 
the vector $f$ has rank exactly equal to $1$. We remark that
the matrices $A_i := L + D_i$ are a natural candidate to build 
a preconditioner using exponential sums. 
\medskip 

\subsubsection{High-dimensional Markov chains}
\label{sec:case-study-markov} Our 
second test case arises from the description of a 
Markov chain. The case study we describe is often
found when dealing with the evaluation of performance
and reliability measures of complex systems, for which 
a high-dimensional state-space naturally appears. Consider 
a set of $d$ systems that evolve stochastically as a 
continuous time Markov 
chain, each of them endowed 
with a state-space $\mathcal S_i$, with 
$|\mathcal S_i| = n$. Even though 
the combined state space would be $\mathcal S := \prod_{i = 1}^d \mathcal S_i$, which has cardinality $n^d$, this high-dimensional Markov chain is relatively easy to analyze because every system evolves independently of the others. 

We now modify the Markov chain allowing some state transitions inside $\mathcal S$ that involve more than 
one system (called \emph{synchronizations}). This 
situation may arise for instance when analyzing computer
networks, where failure of one server may impact one 
or more other services. With this modification, 
the systems cannot be analyzed independently anymore, 
and the problem is truly high-dimensional. The computation
of the steady-state probabilities can be recast to 
solving a linear system of the form 
\[
  (Q + W - D) \pi = e, \qquad 
  Q = \bigoplus_{i = 1}^d Q_i,
\]
where $Q_i$ encodes the transition rates of the systems 
when viewed independently, $W$ adds the synchronization
transitions, and $D$ is a diagonal matrix to ensure 
that the row-sum is zero. The vector $\pi$ contains
the steady-state probabilities. 

This kind of system has been previously analyzed in \cite{kaes,perf1}. We refer the interested reader to these 
works and the references therein for further
details on the model. In this work, we assume that we 
have a family of $d$ systems with the following 
interaction topology:
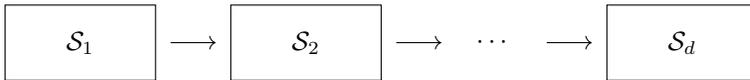
\begin{center}
\medskip 
\begin{tikzpicture}
    \draw (0,0) rectangle (2, 1);
    \node at (1, 0.5) {$\mathcal S_1$};
    \draw[->] (2.2, .5) -- (2.8, .5);
    \draw (3,0) rectangle (5, 1);
    \node at (4, 0.5) {$\mathcal S_2$};
    \draw[->] (5.2, .5) -- (5.8, .5);
    \node at (6.5, .5) {$\cdots$};
    \draw[->] (7.2, .5) -- (7.8, .5);
    \draw (8,0) rectangle (10, 1);
    \node at (9, 0.5) {$\mathcal S_d$};
\end{tikzpicture}
\medskip 
\end{center}
We assume that when particular transitions in system $\mathcal S_i$
are triggered, they change the state in system $\mathcal S_{i+1}$, for all $i < d$. As mentioned above, these particular transitions are called \emph{synchronizations}. Note that this fits well 
with the underlying topology of indices in tensor-trains, 
and often allows to represent the steady-state vector 
in this low-rank format efficiently. 
The transition rates are chosen as 
follows:
\begin{itemize}
    \item Each system behaves as a random walk, 
      with transition rates $\eta_k$ and $\mu_k$ to move forward and 
      backward from state $k$ 
      chosen with a random uniform distribution from 
      $[1, 2]$. All transition rates are chosen 
      independently (that is, the systems 
      are not equidistributed).     
    \item Systems $i$ and $i+1$ have a synchronized 
      transition such that when both systems are
      in state $n-1$, they move together to state $n$ (in the model, this represents the failure of both systems at once). The rate of ``joint failure'' 
      is equal to $0.1$ in our model.
\end{itemize}
From the linear algebra point of view, this means that 
the matrices $Q_i$ are all tridiagonal, and $W$ is the sum of matrices
obtained by the Kronecker product of $d - 2$ identity 
matrices (corresponding to the systems not impacted by
the failure) and $2$ matrices with only one non-zero entry. 

\begin{remark}
    The sparse 
    structure of the matrices could be exploited 
    for both case studies in sections~\ref{sec:case-study-convection-diffusion} and \ref{sec:case-study-markov} to
    accelerate the matvec operations. For the sake of 
    simplicity, generality, 
    and readability of the code we avoided 
    doing so, but we expect that this could be a 
    further speed-up to our experiments. 
\end{remark}

% \rr{@Leo: inserire descrizione PDE problem}

\subsection{Unpreconditioned GMRES}
\label{sec:numexp-unpreconditioned}

In this section, we analyze the performances of TT-sGMRES without preconditioning, applied to the two nonsymmetric problems described above: the convection-diffusion case study and the Markov chain one. In these problems, 
the condition number depends polynomially on $n$, and therefore we only 
considerd small values of $n$, and test the scaling with the number of 
dimensions. 

\subsubsection{Loss in accuracy of vanilla TT-sGMRES}

The first experiment has the aim of showing that the ``vanilla'' 
TT-sGMRES presented Algorithm~\ref{alg:sTT-GMRES} has accuracy problems in the 
reconstruction of the solution, whereas this is not the case in 
the ``enhanced'' TT-sGMRES that we presented in Algorithm~\ref{alg:enhanced sTT-GMRES}. In fact, since the matrix $W_k$ obtained by
running Arnoldi with partial reorthogonalization 
becomes increasingly poorly conditioned, we expect to find large 
cancellations when reconstructing the final solution. This leads to 
poor accuracy if successive relative truncations are performed while 
computing the sum, which are instead avoided when approximating the 
sum all at once with the STTA scheme of section~\ref{Randomized approximation of sums in TT-format}. 

\begin{figure}
\centering
\begin{tikzpicture}
\begin{semilogyaxis}[
    x label style={at={(axis description cs:0.5,0)},anchor=north},
    xlabel={{Iteration}},
    xlabel style = {yshift=-10},
    ylabel={{$\frac{\|Ax_k-b\|_2}{\|b\|_2}$}},
    legend pos=north west,
    ymajorgrids=true,
    grid style=dashed,
    width=.55\linewidth,
    legend pos = north east,
    mark size = 1.2pt,
    restrict x to domain=0:500
]

\addplot [color = red, style = thick,mark = square*,  mark size=0.1pt] table [col sep=tab, x index = 0, y index = 1,x filter/.code={\pgfmathparse{mod(\coordindex,50)==0 ? x : inf}}] {ex_0_res_vs_sres.dat};
\addplot [color = blue, style = thick,mark = triangle*,  mark size=0.1pt] table [col sep=tab, x index = 0, y index = 2,x filter/.code={\pgfmathparse{mod(\coordindex,50)==0 ? x : inf}}] {ex_0_res_vs_sres.dat};

\legend{\scriptsize{Vanilla - Alg.~\ref{alg:sTT-GMRES}}, \scriptsize{Enhanced - Alg.~\ref{alg:enhanced sTT-GMRES}}}
\end{semilogyaxis}
\end{tikzpicture}
 \caption{Actual residuals of the vanilla and enhanced TT-sGMRES algorithms computed after each iteration for the PDE problem in section \ref{sec:case-study-convection-diffusion} with $d = 4$ and $n=34$.} 
    \label{fig:vanilla-vs-enhanced}
\end{figure}
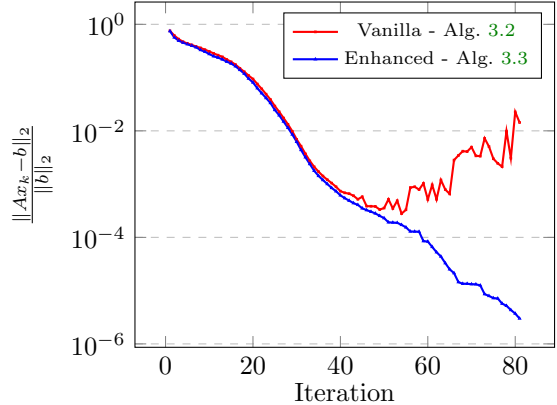

For this case we set $\ell = 1$, and run the vanilla and enhanced version 
of TT-sGMRES on the same problem with $n = 34$ and $d = 4$, for $80$ iterations. The two algorithms are exactly the same with the only exception 
of the final reconstruction described in line \ref{alg:final_line22} of Algorithm \ref{alg:sTT-GMRES}. We then show the value of the residual 
(recomputed exactly) at each iteration, and report it for both 
schemes in Figure~\ref{fig:vanilla-vs-enhanced}. While the enhanced 
version shows a nice convergence plot, the vanilla one has a semiconvergent
behavior, and starting from iteration $40$ the cancellation errors completely
dominate with respect to the achieved accuracy. 
%------------------------- EX vanilla-vs-enhanced

Since the one depicted in Figure~\ref{fig:vanilla-vs-enhanced} is a common behavior of the vanilla TT-sGMRES in the following we focus only on Algorithm~\ref{alg:enhanced sTT-GMRES}. 

%--------------------------- EX1
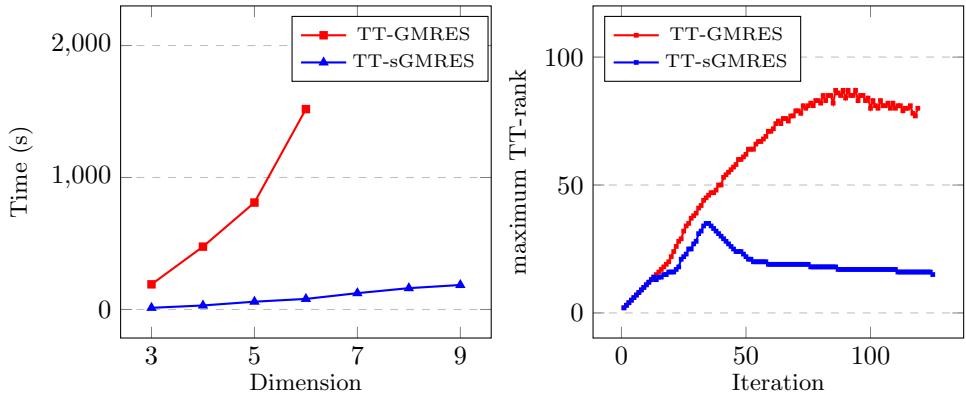
\begin{figure}
\centering
\begin{tikzpicture}
\begin{axis}[
    x label style={at={(axis description cs:0.5,0)},anchor=north},
    xlabel={\small{Dimension}},
    xlabel style = {yshift=-10},
    ylabel={\small{Time (s)}},
    ymajorgrids=true,
    grid style=dashed,
    width=.5\linewidth,
    height = 6cm, 
    ymax = 2300, 
    legend style={scale=0.8, xshift = 4}, 
    legend pos = north east,
    mark size = 1.2pt,
    restrict x to domain=0:125,
    xtick = {3,5,7,9}
]

\addplot [color = red, style = thick,mark = square*,  mark size=1.2pt] table [col sep=tab, x index = 0, y index = 1,x filter/.code={\pgfmathparse{mod(\coordindex,50)==0 ? x : inf}}] {ex_1_times_tt.dat};
\addplot [color = blue, style = thick,mark = triangle*,  mark size=1.5pt] table [col sep=tab, x index = 0, y index = 1,x filter/.code={\pgfmathparse{mod(\coordindex,50)==0 ? x : inf}}] {ex_1_times_stt.dat};

\legend{\scriptsize{TT-GMRES}, \scriptsize{TT-sGMRES}}
\end{axis}
\end{tikzpicture}~\begin{tikzpicture}
%%%%%%%%%%%%%%% EX2
\begin{axis}[
    x label style={at={(axis description cs:0.5,0)},anchor=north},
    xlabel={\small{Iteration}},
    xlabel style = {yshift=-10},
    ylabel={\small{maximum TT-rank}},
    legend pos=north west,
    ymajorgrids=true,
    grid style=dashed,
    width=.5\linewidth,
    ymax = 120, 
    height = 6cm, 
    legend pos = north west,
    mark size = 1.2pt,
    restrict x to domain=0:125
]

\addplot [color = red, style = thick,mark = square*,  mark size=0.5pt] table [col sep=tab, x index = 0, y index = 1,x filter/.code={\pgfmathparse{mod(\coordindex,50)==0 ? x : inf}}] {ex_2_tt_ranks.dat};
\addplot [color = blue, style = thick,mark = square*,  mark size=0.5pt] table [col sep=tab, x index = 0, y index = 1,x filter/.code={\pgfmathparse{mod(\coordindex,50)==0 ? x : inf}}] {ex_2_stt_ranks.dat};

\legend{\scriptsize{TT-GMRES}, \scriptsize{TT-sGMRES}}
\end{axis}
\end{tikzpicture}
 \caption{On the left, we report the runtime of the TT-GMRES and TT-sGMRES algorithms on convection-diffusion PDE problems of size $n=64$ across various dimensions $d$ and accuracy $10^{-4}$. On the right, 
 we plot the maximum TT-ranks of the base vectors generated by TT-GMRES and TT-sGMRES with $d=6$, $n=64$ and ${\texttt{tol}=10^{-4}}$. In the
 right experiment, TT-GMRES converged in 1528.22 seconds with respect to the 80.03 seconds of TT-sGMRES.} 
    \label{fig:ex_1_times}
\end{figure}

\subsubsection{TT-GMRES vs TT-sGMRES}
In the second experiment we consider again the PDE problem
from section~\ref{sec:case-study-convection-diffusion}, and we compare
the timings of the enhanced TT-sGMRES with the standard TT-GMRES. The problem
is considered for $d$ ranging from $3$ to $9$, and $n$ fixed to $64$. 
The stopping criterion is ${\texttt{tol}}=10^{-4}$, and we aborted the 
execution if the runtime exceeded one hour. The results are reported in Figure~\ref{fig:ex_1_times} (left). 

In this test, the enhanced TT-sGMRES is faster than TT-GMRES for all dimensions. The speedup arises from two phenomena: we only perform partial 
 reorthogonalization and the TT-ranks remain smaller.
 To better describe the latter phenomenon we provide another plot in Figure \ref{fig:ex_1_times} (right), in which we show the maximum TT-rank of the vectors $v_k$ generated by the two algorithms for $d=5$ (for other dimensions we obtained analog results). We can see that TT-GMRES operates with higher TT-ranks with respect to the enhanced TT-sGMRES. On one side, higher TT-ranks lead to more expensive arithmetic operations, and on the other side the fact that TT-GMRES performs full orthogonalization increases the number of dot products; the enhanced TT-sGMRES, instead, 
 only requires a constant number of these dot products per 
 iteration. We also observe that the enhanced TT-sGMRES requires a few more iterations to converge than TT-GMRES, mostly because the sketched tolerance is set to $0.3\cdot{\texttt{tol}}$ in order to accommodate with 
 the estimation error for the residual; however, in most cases the 
 sketched residual is a relatively tight estimate, so we end with 
 a slightly more accurate solution with TT-sGMRES than with 
 TT-GMRES.

\subsubsection{Gap between sketched and actual residual}
In the previous examples we have set the tolerance for the stopping criterion 
in TT-sGMRES slightly smaller than the one for TT-GMRES. This is because 
the stopping criterion for the former relies on the sketched relative
residual $\| S(Ax_k - b) \| / \| Sb \|$, which is a good estimate of the 
true residual up to a small constant (with high probability). 

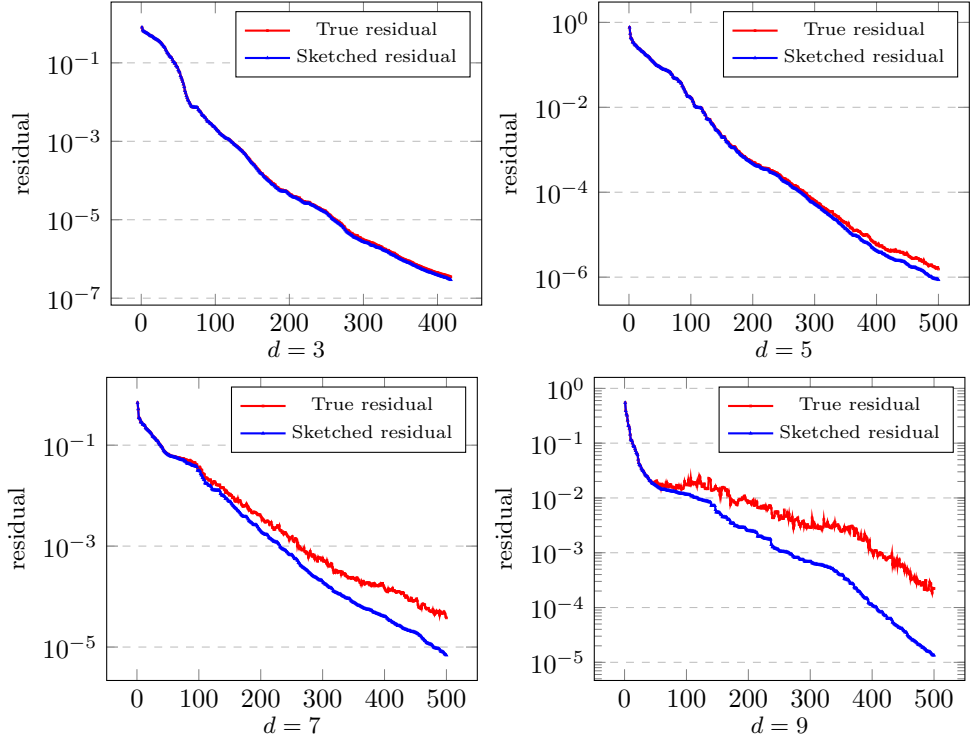
\begin{figure}
\centering
\begin{tikzpicture}
\begin{semilogyaxis}[
    x label style={at={(axis description cs:0.5,0)},anchor=north},
    xlabel={\small{$d=3$}},
    xlabel style = {yshift=-10},
    ylabel={\small{residual}},
    legend pos=north west,
    ymajorgrids=true,
    grid style=dashed,
    width=.50\linewidth,
    legend pos = north east,
    mark size = 1.2pt,
    restrict x to domain=0:500
]

\addplot [color = red, style = thick,mark = square*,  mark size=0.1pt] table [col sep=tab, x index = 0, y index = 1,x filter/.code={\pgfmathparse{mod(\coordindex,50)==0 ? x : inf}}] {ex_3_res_vs_sres_dim_3.dat};
\addplot [color = blue, style = thick,mark = triangle*,  mark size=0.1pt] table [col sep=tab, x index = 0, y index = 2,x filter/.code={\pgfmathparse{mod(\coordindex,50)==0 ? x : inf}}] {ex_3_res_vs_sres_dim_3.dat};

\legend{\scriptsize{True residual}, \scriptsize{Sketched residual}}
\end{semilogyaxis}
\end{tikzpicture}~\begin{tikzpicture}
\begin{semilogyaxis}[
    x label style={at={(axis description cs:0.5,0)},anchor=north},
    xlabel={\small{$d=5$}},
    xlabel style = {yshift=-10},
    ylabel={\small{residual}},
    legend pos=north west,
    ymajorgrids=true,
    grid style=dashed,
    width=.5\linewidth,
    legend pos = north east,
    mark size = 1.2pt,
    restrict x to domain=0:500
]

\addplot [color = red, style = thick,mark = square*,  mark size=0.1pt] table [col sep=tab, x index = 0, y index = 1,x filter/.code={\pgfmathparse{mod(\coordindex,50)==0 ? x : inf}}] {ex_3_res_vs_sres_dim_5.dat};
\addplot [color = blue, style = thick,mark = triangle*,  mark size=0.1pt] table [col sep=tab, x index = 0, y index = 2,x filter/.code={\pgfmathparse{mod(\coordindex,50)==0 ? x : inf}}] {ex_3_res_vs_sres_dim_5.dat};

\legend{\scriptsize{True residual}, \scriptsize{Sketched residual}}
\end{semilogyaxis}
\end{tikzpicture}\\

\begin{tikzpicture}
\begin{semilogyaxis}[
    x label style={at={(axis description cs:0.5,0)},anchor=north},
    xlabel={\small{$d=7$}},
    xlabel style = {yshift=-10},
    ylabel={\small{residual}},
    legend pos=north west,
    ymajorgrids=true,
    grid style=dashed,
    width=.5\linewidth,
    legend pos = north east,
    mark size = 1.2pt,
    restrict x to domain=0:500
]

\addplot [color = red, style = thick,mark = square*,  mark size=0.1pt] table [col sep=tab, x index = 0, y index = 1,x filter/.code={\pgfmathparse{mod(\coordindex,50)==0 ? x : inf}}] {ex_3_res_vs_sres_dim_7.dat};
\addplot [color = blue, style = thick,mark = triangle*,  mark size=0.1pt] table [col sep=tab, x index = 0, y index = 2,x filter/.code={\pgfmathparse{mod(\coordindex,50)==0 ? x : inf}}] {ex_3_res_vs_sres_dim_7.dat};

\legend{\scriptsize{True residual}, \scriptsize{Sketched residual}}
\end{semilogyaxis}
\end{tikzpicture}~\begin{tikzpicture}
\begin{semilogyaxis}[
    x label style={at={(axis description cs:0.5,0)},anchor=north},
    xlabel={\small{$d=9$}},
    xlabel style = {yshift=-10},
    ylabel={\small{residual}},
    legend pos=north west,
    ymajorgrids=true,
    grid style=dashed,
    width=.5\linewidth,
    legend pos = north east,
    mark size = 1.2pt,
    restrict x to domain=0:500
]

\addplot [color = red, style = thick,mark = square*,  mark size=0.1pt] table [col sep=tab, x index = 0, y index = 1,x filter/.code={\pgfmathparse{mod(\coordindex,50)==0 ? x : inf}}] {ex_3_res_vs_sres_dim_9.dat};
\addplot [color = blue, style = thick,mark = triangle*,  mark size=0.1pt] table [col sep=tab, x index = 0, y index = 2,x filter/.code={\pgfmathparse{mod(\coordindex,50)==0 ? x : inf}}] {ex_3_res_vs_sres_dim_9.dat};

\legend{\scriptsize{True residual}, \scriptsize{Sketched residual}}
\end{semilogyaxis}
\end{tikzpicture}
 \caption{The above plots report the difference between the 
 sketched residual and the true residual, for different values of 
 $d$.} 
    \label{fig:ex_3_res_vs_sres}
\end{figure}

In this experiment, we show the distance between the sketched and the true residuals, for various dimensions $d=3,5,7,9$. The results along all the iterations for the 
PDE problem with $n = 64$ are reported in Figure~\ref{fig:ex_3_res_vs_sres}.
The maximum number of iterations is set to $500$, and the number of rows 
of $S$ to $1000$, so at the end of the algorithm the dimension of the 
sketched space is about twice as the dimension of the subspace where the 
residual lives. The tolerance was set to ${\texttt{tol}}=10^{-6}$.

The plots show that the gaps are higher for higher values of $d$. One of the causes is that the embedding power of the Khatri-Rao embeddings depends on the dimension $d$, the other, and most impactful, is that STTA recovers an approximate low-rank approximation up to some constants depending exponentially on $d$.

In this experiment and in the tests that we have run, this gap has always been less than $10$; however, for higher dimensions, this gap could become significant, because of the loss of accuracy of the STTA approximation. It is possible to compensate this effect and reduce the STTA constants by increasing the parameter $\ell_\mu$ in the generation of the sketchings phase. For further details, see \cite{kressner2023streaming}.

%--------------------------------- EX3

\subsubsection{Markov case study without preconditioning}

\begin{figure}
\centering
\begin{tikzpicture}
\begin{axis}[
    x label style={at={(axis description cs:0.5,0)},anchor=north},
    xlabel={\small{Dimension}},
    xlabel style = {yshift=-10},
    ylabel={\small{Time (s)}},
    legend pos=north west,
    ymajorgrids=true,
    grid style=dashed,
    width=.48\linewidth,
    legend pos = north west,
    mark size = 1.2pt,
    restrict x to domain=0:125
]

\addplot [color = red, style = thick,mark = square*] table [col sep=tab, x index = 0, y index = 2,x filter/.code={\pgfmathparse{mod(\coordindex,50)==0 ? x : inf}}] {ex_5_time.dat};
\addplot [color = blue, style = thick,mark = triangle*] table [col sep=tab, x index = 0, y index = 1,x filter/.code={\pgfmathparse{mod(\coordindex,50)==0 ? x : inf}}] {ex_5_time.dat};

\legend{\scriptsize{TT-GMRES}, \scriptsize{TT-sGMRES}}
\end{axis}
\end{tikzpicture}~\begin{tikzpicture}
\begin{axis}[
    x label style={at={(axis description cs:0.5,0)},anchor=north},
    xlabel={\small{Iteration}},
    xlabel style = {yshift=-10},
    ylabel={\small{maximum TT-rank}},
    ymajorgrids=true,
    grid style=dashed,
    width=.48\linewidth,
    legend pos = north east,
    mark size = 1.2pt,
    restrict x to domain=0:125,
    ymax = 200
]

\addplot [color = red, style = thick,mark = square*,  mark size=0.5pt] table [col sep=tab, x index = 0, y index = 1,x filter/.code={\pgfmathparse{mod(\coordindex,50)==0 ? x : inf}}] {ex_5_tt_ranks_dim_5.dat};
\addplot [color = blue, style = thick,mark = square*,  mark size=0.5pt] table [col sep=tab, x index = 0, y index = 1,x filter/.code={\pgfmathparse{mod(\coordindex,50)==0 ? x : inf}}] {ex_5_stt_ranks_dim_5.dat};

\legend{\scriptsize{TT-GMRES}, \scriptsize{TT-sGMRES}}
\end{axis}
\end{tikzpicture}
 \caption{On the left, the comparison between running 
 TT-GMRES and TT-sGMRES for the Markov test case, with different values 
 of $d$ and $n = 64$. On the right, the behavior of ranks of the 
 basis vectors during the iterations, in the case $d = 5$.} 
    \label{fig:ex_4_ranks}
\end{figure}
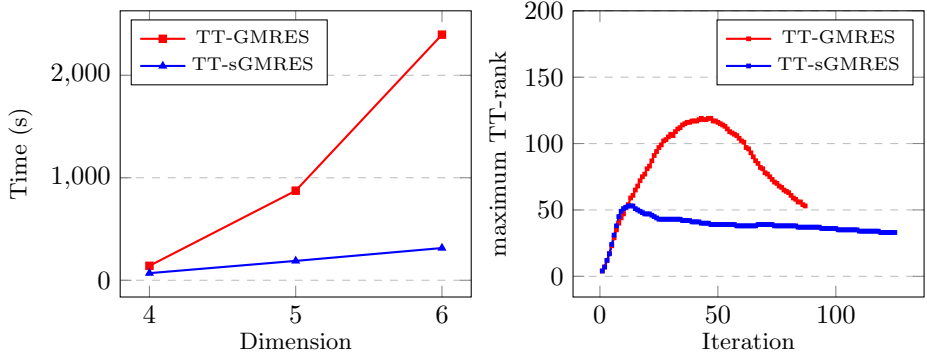

We have replicated the experiments for the PDE problems on the Markov 
case study, which led to a similar behavior. We report in this 
section the timings for running TT-GMRES and TT-sGMRES, which are 
plotted in Figure~\ref{fig:ex_4_ranks}, on the left. We can see that, 
as in the PDE case study, the proposed algorithm can deal with the 
increasing dimensionality without a significant increase in 
computational times (with respect to TT-GMRES). 

On the right, 
in the same Figure, the ranks throughout the iterations are reported. 
In contrast with the PDE example, the rank of the operator describing the 
Markov chain grows with $d$ (linearly), and therefore the problem
becomes increasingly challenging for high dimensions. 

We remark, however, that without preconditioning the performances of the algorithm are still far from those of AMEn. Therefore in the next section, we focus on the preconditioned case.

%--------------------------------- EX4

% %--------------------------------- EX6

% \begin{figure}
% \centering
% \begin{tikzpicture}
% \begin{axis}[
%     x label style={at={(axis description cs:0.5,0)},anchor=north},
%     xlabel={\small{Dimension}},
%     ylabel={\small{Rank}},
%     legend pos=north west,
%     ymajorgrids=true,
%     grid style=dashed,
%     width=.6\linewidth,
%     legend pos = north east,
%     mark size = 1.2pt,
%     restrict x to domain=0:400
% ]

% \addplot [color = red, style = thick,mark = square*,  mark size=0.5pt] table [col sep=tab, x index = 0, y index = 1,x filter/.code={\pgfmathparse{mod(\coordindex,50)==0 ? x : inf}}] {ex_6_tt_ranks.dat};
% \addplot [color = blue, style = thick,mark = square*,  mark size=0.5pt] table [col sep=tab, x index = 0, y index = 1,x filter/.code={\pgfmathparse{mod(\coordindex,50)==0 ? x : inf}}] {ex_6_stt_ranks.dat};

% \legend{\scriptsize{$TT$}, \scriptsize{$sTT$}}
% \end{axis}
% \end{tikzpicture}
%  \caption{EX6: \rr{A:questo plot potremmo toglierlo, a meno che non vogliamo più esperimenti di questo tipo}Ranks comparison between TT-gmres and STT gmres on Markov problem, dimension 6 and with $n=16$, tt gmres spent 893 seconds with respect to the 313 of sTTGMRES (which also run to higher accuracy $0.3*\mathrm{tol}$ which is almost preserved)} 
%     \label{fig:ex_6_ranks}
% \end{figure}

\subsection{Numerical tests with preconditioning}
\label{sec:numexp-preconditioned}
In this section, we reconsider the case studies 
presented above, and include an option to precondition
the TT-sGMRES iteration. In both cases, this is necessary 
when the dimensions $n_i$ become large, because the 
condition number grows polynomially in $n$. We will use 
exponential sums to build preconditioners for all examples 
for simplicities, but we do not expect major differences in 
case other preconditioners are used. 

\subsubsection{Convection-diffusion}
For the convection-diffusion problem in the case $d = 5$, we employed 
an exponential sum preconditioner with 
\[
  \mathcal P = \sum_{0}^{\zeta} \alpha_j \bigotimes_{i = 1}^d e^{-\beta_j A_i}, 
\]
as detailed in section~\ref{sec: preconditioning}. We have 
chosen to take $\zeta = 17$. In addition, we have tested 
Algorithm~\ref{alg:Preconditioned_sTT-GMRES} with 
different values of \texttt{maxrank}. As a rule of 
thumb, we expect smaller values of \texttt{maxrank}
to yields faster iterations, but slower convergence, or even stagnation. On the other hand, higher values of \texttt{maxrank} will be closer to the GMRES
iteration without rounding and usually yield a better 
convergence, but with a much higher computational 
cost per iteration. 

For this example, we have tested 
\texttt{maxrank} $= \infty$ and 
\texttt{maxrank} $= 30$; in addition, we have compared 
the performances with the AMEn solver in the 
TT-Toolbox (with default parameters, and a maximum 
number of sweeps set to $200$ in order to achieve the target tolerance). The target tolerance was set to 
$10^{-8}$, and as usual we reduced it by a factor 
$10$ in TT-sPGMRES, to account for the 
constant in the estimation of the residual by 
sketching. 

\begin{figure}
    \centering
    \begin{tikzpicture}
        \begin{loglogaxis}[legend pos = north west, 
          width = 10cm, height = 6cm,
          xtick = {128, 256, 512, 1024}, 
          ylabel = {Time (s)}, 
          xlabel = {$n_i$}, 
          xticklabel={\pgfmathparse{exp(\tick)}\pgfmathprintnumber[precision=0]{\pgfmathresult}},
          x label style={at={(axis description cs:0.5,0)},anchor=north},
          ymajorgrids=true,
          grid style=dashed,
          mark size = 1.2pt,
]

            \addplot[color = red, style = thick, mark = square*] 
               table[x index = 0, y index = 4] {ex_precond_1.dat};
            \addplot[color = blue, style = thick, mark = triangle*] 
               table[x index = 0, y index = 7] {ex_precond_1.dat};
            \addplot[color = black, style = thick, mark = diamond*]  table[x index = 0, y index = 10] {ex_precond_1.dat};
            \legend{
                \texttt{maxrank} $= \infty$, 
                \texttt{maxrank} $= 30$,
                AMEn
            };
            
        \end{loglogaxis}
    \end{tikzpicture}
    \caption{Runtime of TT-sPGMRES iteration 
    for the convection-diffusion problem in section~\ref{sec:case-study-convection-diffusion} with variable $n_i$ and $d = 5$; 
    the target tolerance in this example is $10^{-8}$, and 
different values of {\texttt{maxrank}} are used. AMEn is run with standard 
parameters, and is taken from TT-Toolbox \cite{TT-Toolbox}.}
    \label{fig:precond_1}
\end{figure}
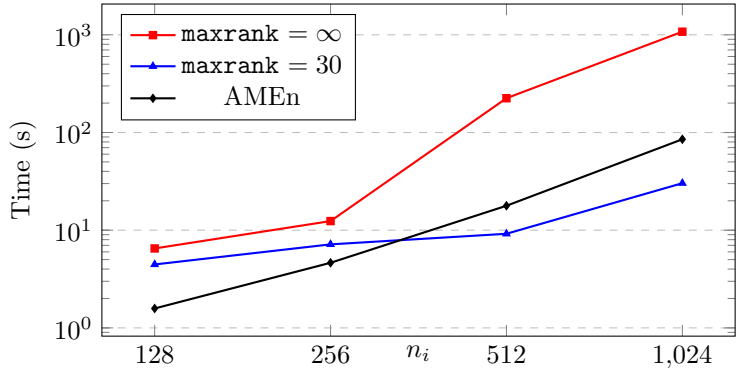

All approaches achieved the required accuracy, and 
the timings for different 
values of $n_i$ are reported in Figure~\ref{fig:precond_1}. 
We see from the results in Figure~\ref{fig:precond_1} that allowing the 
ranks to grow unbounded does not yield optimal performances. With both 
{\texttt{maxrank}} set to $\infty$ and $30$, TT-sPGMRES convergences 
in $4$ iterations to the desired tolerance with this choice of 
preconditioner. Moreover, when choosing 
{\texttt{maxrank}} $= 30$ our algorithm becomes competitive and, for this example, it is faster than AMEn.

Without preconditioning, the ranks stay nicely bounded, but the number 
of iterations is so large that the method cannot be competitive with the 
choices above. With {\texttt{maxrank}} $= \infty$, the iteration 
reaches rank $433$ for $n_i = 1024$, so it is rather memory 
demanding. Hence, this example shows how using a bounded 
rank can be essential when incorporating preconditioning.

\subsubsection{Preconditioning for the Markov test case}

We have run a similar experiment for the test case arising
from Markov chains. In that case, a natural choice is to
consider the infinitesimal generator $Q$ obtained by 
ignoring all interactions between the different systems, and dropping 
the matrix $W$ (following the notation used in section~\ref{sec:case-study-markov}). 

The matrix $Q$ is a Kronecker sum, and therefore its approximate inverse 
can be constructed by exponential sums, exactly as for 
the convection-diffusion test case. For this problem, 
we have selected $\zeta = 33$. We have then run the same 
tests, using systems with a number of states ranging from 
$128$ to $512$, and requiring tolerance $10^{-6}$. This problem 
is more challenging than the PDE case, and we have run our algorithm 
with {\texttt{maxrank}} $\in \{ 50, 80, \infty \}$. As 
for the PDE case, we use $\eta = 0.1$ as a safety factor to make 
sure that the if the sketched residual is below $\eta \cdot \epsilon$
then the true residual is around $\epsilon$ or less. 

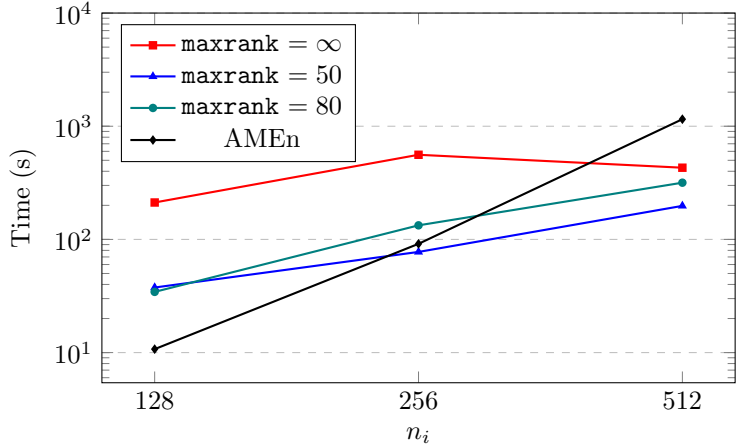
\begin{figure}
    \centering
    \begin{tikzpicture}
        \begin{loglogaxis}[legend pos = north west, 
          width = 10cm, height = 6.5cm,
          xtick = {128, 256, 512, 1024}, 
          ylabel = {Time (s)}, 
          ymax = {10000}, 
          xlabel = {$n_i$}, 
          xticklabel={\pgfmathparse{exp(\tick)}\pgfmathprintnumber[precision=0]{\pgfmathresult}},          
          ymajorgrids=true,
          grid style=dashed,
          mark size = 1.2pt,
]
            \addplot[color = red, style = thick, mark = square*] table[x index = 0, y index = 5] {ex_precond_2.dat};
            \addplot[color = blue, style = thick, mark = triangle*]  table[x index = 0, y index = 9] {ex_precond_2.dat};
            \addplot[color = teal, style = thick, mark = otimes*] table[x index = 0, y index = 13] {ex_precond_2.dat};
            \addplot[color = black, style = thick, mark = diamond*] table[x index = 0, y index = 17] {ex_precond_2.dat};
            \legend{
                \texttt{maxrank} $= \infty$, 
                \texttt{maxrank} $= 50$,
                \texttt{maxrank} $= 80$,
                AMEn
            };
        \end{loglogaxis}
    \end{tikzpicture}
    \caption{Runtime of TT-sPGMRES iteration 
    for the Markov problem in section~\ref{sec:case-study-markov} with variable $n_i$ and $d = 5$; the target tolerance in this example is $10^{-6}$, and 
different values of {\texttt{maxrank}} are used. AMEn is run with standard 
parameters, and is taken from TT-Toolbox \cite{TT-Toolbox}.}
    \label{fig:precond_2}
\end{figure}

When running with {\texttt{maxrank} $= \infty$} we encounter the same behavior 
of the PDE case study of the previous section: 
the rank grows quickly (up to about $220$
in this example), the algorithm is slowed down and can easily encounter memory 
issues. On the other hand, using lower values of 
{\texttt{maxrank}} makes the algorithm competitive with AMEn, and even faster 
for large values of $n_i$, and corresponding badly conditioned problems. In this 
example, {\texttt{maxrank} $= 50$} only manages to reach a true accuracy 
of about $10^{-5}$, whereas {\texttt{maxrank} $= 80$} achieves the target of 
$10^{-6}$.

% \subsection{Comparison with TT-GMRES, AMEn, altro?}

\section{Conclusions}\label{sec: conclusions}
In this work, we presented and analyzed a sketched version of 
TT-GMRES, called TT-sGMRES, a novel algorithm that combines the winning strategies of sketch GMRES and TT-GMRES. Through various methodological refinements, we demonstrated that the introduction of sketching and randomization brings significant benefits, primarily by greatly reducing the cost of orthogonalization and limiting the 
ranks of tensors during the iteration.
Additionally, the approach based on a streamable method allowed us to overcome one of the classic storage problems, namely the allocation of the whole basis. In particular, once the vectors of the 
Krylov bases are computed, they are sketched and then discarded, and this is sufficient to recover the solution upon convergence. 

The experiments conducted validate the effectiveness of the proposed method. Not only did the TT-sGMRES prove to be significantly superior to the classical TT-GMRES, but in many cases, it was also competitive with established solvers such as AMEn.
Another advantage of our method is the possibility of leveraging preconditioners to further improve its performance, making it an extremely promising method for a wide range of applications. 

Although we focused on the TT-format, many of the improvements introduced can be tested and exploited in a broader range of cases where vectors can be compressed in a low-rank format and streamable algorithms for their linear combinations are available. 
For example, this approach could be applied to the Tucker format using the methods in \cite{bucci2023multilinear,bucci2024sequential,sun2020low}, and efforts could be made to extend it to the case of the Tree Tensor Network format.

In conclusion, TT-sGMRES represents a significant advancement in the state of the art, offering an efficient and scalable scheme for solving high-dimensional linear systems.

\section*{Acknowledgments}
The authors are members of the INdAM Research
Group GNCS that partially supported this work through the funded project GNCS2024 with reference number CUP\_E53C23001670001.

The work of the authors
was partially supported also by the European Union - NextGenerationEU under the National Recovery and Resilience Plan (PNRR) - Mission 4 Education and research
- Component 2 From research to business - Investment 1.1 Notice Prin 2022 - DD N. 104 of 2/2/2022,
entitled “Low-rank Structures and Numerical Methods in Matrix and Tensor Computations and their
Application”, code 20227PCCKZ – CUP J53D23003620006.

\bibliographystyle{siamplain}
\bibliography{references}
\end{document}